\documentclass[12pt,reqno,a4paper]{amsart}
\usepackage[left=1.5cm,right=1.5cm,top=3cm,bottom=3cm]{geometry}

\usepackage{cite}
\usepackage{amssymb,amsfonts,amsmath}
\usepackage{mathtools} 
\usepackage{color}
\usepackage{graphicx,float}

\usepackage{booktabs}
\usepackage{dcolumn,cancel}
\newcolumntype{d}[1]{D{.}{.}{#1}}  
\usepackage[normalem]{ulem}

\newtheorem{theorem}{Theorem}
\newtheorem{lemma}{Lemma}
\newtheorem{assum}{Assumption}

\newtheorem{remark}{Remark}

\newcommand{\eps}{\epsilon}

\newcommand{\be}{\boldsymbol{e}}

\newcommand{\bG}{\boldsymbol{G}}
\newcommand{\bI}{\boldsymbol{I}}

\newcommand{\bk}{{\boldsymbol{k}}}
\newcommand{\ba}{{\boldsymbol{a}}}

\newcommand{\bmu}{{\boldsymbol{\mu}}}
\newcommand{\brho}{{\boldsymbol{\rho}}}
\newcommand{\bchi}{{\boldsymbol{\chi}}}
\newcommand{\boldeta}{{\boldsymbol{\eta}}}

\newcommand{\bx}{\boldsymbol{x}}

\newcommand{\by}{\boldsymbol{y}}
\newcommand{\bg}{\boldsymbol{g}}
\newcommand{\br}{\boldsymbol{r}}
\newcommand{\bth}{\boldsymbol{\theta}}
\newcommand{\bzero}{\boldsymbol{0}}
\newcommand{\bgamma}{\boldsymbol{\gamma}}
\newcommand{\bnu}{\boldsymbol{\nu}}

\newcommand{\hh}{\hspace*{0.7pt}}
\newcommand{\nes}{\hspace*{-0.2pt}}
\newcommand{\ww}{w}
\newcommand{\tw}{\tilde{w}}
\newcommand{\pho}{\tilde{\phi}_n^{\circ}}
\newcommand{\pha}{\tilde{\phi}_n^{\ba}}
\newcommand{\vpha}{\tilde{\varphi}_n^{\ba}}
\newcommand{\vphah}{\tilde{\varphi}_{\hat{n}}^{\ba}}
\newcommand{\phaq}{\tilde{\phi}_{nq}^{\ba}}
\newcommand{\vphaq}{\tilde{\varphi}_{nq}^{\ba}}
\newcommand{\vphao}{\tilde{\varphi}_{n1}^{\ba}}

\newcommand{\vphap}{\tilde{\varphi}_{np}^{\ba}}
\newcommand{\vphar}{\tilde{\varphi}_{nr}^{\ba}}
\newcommand{\vphas}{\tilde{\varphi}_{ns}^{\ba}}

\begin{document}

\title[Dynamic homogenization at finite wavelengths and frequencies]{A rational framework for dynamic homogenization at finite wavelengths and frequencies}

\author{Bojan B. Guzina$^{1}$, Shixu Meng$^{2}$ and Othman Oudghiri-Idrissi$^{1}$}

\address{$^{1}$Department of Civil, Environmental, and Geo- Engineering, University of Minnesota, Twin Cities\\
$^{2}$Department of Mathematics, University of Michigan, Ann Arbor}




\begin{abstract}
In this study, we establish an inclusive paradigm for the homogenization of scalar wave motion in periodic media (including the source term) at finite frequencies and wavenumbers spanning the first Brillouin zone. We take the eigenvalue problem for the unit cell of periodicity as a point of departure, and we consider the projection of germane Bloch wave function onto a suitable eigenfunction as descriptor of effective wave motion. For generality the finite wavenumber, finite frequency (FW-FF) homogenization is pursued in~$\mathbb{R}^d$ via second-order asymptotic expansion about the apexes of ``wavenumber quadrants'' comprising the first Brillouin zone, at frequencies near given (acoustic or optical) dispersion branch. We also consider the junctures of dispersion branches and ``dense'' clusters thereof, where the asymptotic analysis reveals several distinct regimes driven by the parity and symmetries of the germane eigenfunction basis. In the case of junctures, one of these asymptotic regimes is shown to describe the so-called Dirac points, that are relevant to the phenomenon of topological insulation. On the other hand, the effective model for nearby solution branches is found to invariably entail a Dirac-like system of equations that describes the interacting dispersion surfaces as ``blunted cones''. For all cases considered, the effective description turns out to admit the same general framework, with differences largely being limited to (i) the eigenfunction basis, (ii)~the reference cell of medium periodicity, and (iii)~the wavenumber-frequency scaling law underpinning the asymptotic expansion. We illustrate the analytical developments by several examples, including the Green's function near the edge of a band gap and clusters of nearby dispersion surfaces.  
\end{abstract}

\maketitle

\section{Introduction} \label{Introduction}

\noindent Effective i.e. ``macroscopic'' description of wave motion in periodic media is a subject of mounting interest in science and engineering for its potential to reduce the computational effort and help physical intuition when tackling problems related to e.g. cloaking, sub-wavelength imaging, and vibration control~\cite{Capo2009,Bava2013}. Naturally, such an idea raises the fundamental question of the field equation governing the effective quantities. One keen approach to the macroscopic wave description in periodic media is the Willis' concept of effective constitutive relationships~\cite{MW07,MBW06,NHA15,NSK12,W11,MG17}.  Another school of thought is the two-scale framework~\cite{PBL78} that has been pursued in the mathematics literature for a long time.

The effective description of wave motion in the long-wavelength, low-frequency \mbox{(LW-LF)} regime is nowadays largely understood. Recent efforts on this front are focused on enriched dispersive models that transcend the usual quasi-static approximation. In this vein, higher-order asymptotic models of the LW-LF motion are often sought via either the Floquet-Bloch analysis or the two-scale homogenization approach~\cite{ABV16,DLB14,CF01,MG17,SS91} to describe incipient wave dispersion within the first, i.e. acoustic pass band. 

Unfortunately some the more intriguing wave phenomena in periodic media, such those of negative refraction~\cite{W16} and topologically protected states~\cite{TWZ17}, occur at finite frequencies that are beyond the reach of the LW-LF asymptotic approximations~\cite{ABV16,CF01,MG17,SS91}. Further, several authors have argued from the platform of Willis' dynamic homogenization that the customary concept of ``mean motion'' (computed by averaging over the unit cell) is not suitable as a macroscopic descriptor at finite frequencies~\cite{NHA15,SN14,NHA16a}. This in particular motivates studies in the low-wavenumber, finite-frequency (LW-FF) regime with the goal of unmasking the effective wave motion and establishing its governing equation. From a broader perspective, early ideas of homogenization at finite frequencies can be found in the asymptotic analysis~\cite{PBL78} of the wave equation with initial data. The leading-order, two-scale homogenization in the LW-FF regime has been particularly considered in more recent studies~\cite{BS06,CKP10,DBDF02} dealing with time-harmonic waves, which exposed the internal resonances of a periodic system as the lynchpin of macroscopic description. To shed further light on the problem, we adopt the framework of plane-wave expansion and pursue \emph{second-order homogenization} of the time-harmonic wave equation in the LW-FF regime via a generalized definition of effective wave motion, given by the projection onto a specific (highly oscillating) Bloch eigenfunction. Such mapping in particular separates the macroscopic and microscopic field fluctuations at finite frequencies, allowing for a systematic asymptotic formulation of the enriched governing equation in~$\mathbb{R}^d$. The analysis demonstrates that the leading correction of the effective differential operator is~$O(\eps^2)$, $\eps$ being the perturbation parameter. In contrast, the effective source term is found to undergo an~$O(\epsilon)$ correction, which is directly relevant to the development of homogenized Green's functions for periodic media. Inspired by the idea of multi-cell homogenization~\cite{Vasi2006,Gone2010,CKP10}, we next extend our analysis to the finite-wavenumber, finite-frequency (FW-FF) regime via expansion about the apexes of ``wavenumber quadrants'' comprising the first Brillouin zone. We also consider (to the leading order) the \emph{junctures of dispersion branches} with arbitrary multiplicity -- relevant to the phenomenon of topologically protected edge states~\cite{TWZ17}, where the effective description is found to depend on the parity and symmetries of the eigenfunction basis affiliated with a repeated eigenvalue. For completeness the latter analysis is generalized to describe the \emph{clusters of nearby solution branches}, including part junctures, where the governing model is found to entail a system of Dirac-like equations -- which paints the local geometry of nearby dispersion surfaces as ``blunted cones''.

All effective descriptions in this work are shown to admit the same general framework, with differences being limited to (a) the Bloch eigenfunction basis, (b) the reference cell of medium periodicity, and (c) the wavenumber-frequency scaling law behind the asymptotic expansion. Consistent with the recent LW-LF analysis~\cite{MG17} we find that the effective FW-FF source terms are, at higher orders of approximation, endowed with respective \emph{correction factors} -- polynomials in frequency and wavenumber -- that account for the interplay between the source density and periodic medium. As an illustration, we apply our analysis toward exposing the effective wave motion in polyatomic chains and chessboard-like continua. This includes an improved asymptotic model for the Green's function near the edge of a band gap~\cite{Doss2008 ,Cra11} -- that highlights the importance of source correction, and effective approximation of closely spaced dispersion surfaces. 

From a broader vantage point this work can be viewed as an extension of the $\bk\!\cdot\!\boldsymbol{p}$ method, e.g.~\cite{Card2005}, that (i) considers the wave motion with a source term; (ii) exposes the leading-order corrections of both effective differential operator and effective source density, and (iii) sheds further light on the interaction of nearby band structures. As to the latter problem, we likewise refer the reader to the theory of invariants, e.g.~\cite{Bir1974}, that makes use of symmetry groups toward a qualitative description of merging or nearby dispersion surfaces.

\section{Preliminaries} \label{prelim}

\noindent With reference to a Cartesian coordinate system endowed with an orthonormal vector basis~$\be_j$ ($j\!=\!\overline{1,d}$), consider the time-harmonic wave equation
\begin{equation}\label{PDE}
-\omega^2\rho(\bx)\hh u - \nabla\!\cdot\!\big(G(\bx)\nabla u\big) ~=~ f(\bx) \quad\text{in~~}\mathbb{R}^d  
\end{equation}
at frequency~$\omega$, where~$G$ and~$\rho$ are $Y$-periodic; 
\begin{equation}\label{Y}
Y \,=\, \{\bx\!:0<\bx\!\cdot\!\be_j < \ell_j; \; j=\overline{1,d}\}, \quad  |Y|=1
\end{equation}
is the unit cell of periodicity, and~$f$ is the source term. In what follows, $G$ and~$\rho$ are assumed to be real-valued $L^\infty(Y)$ functions bounded away from zero. When~$d\!=\!2$, one may interpret~(\ref{PDE}) in the context of anti-plane shear (elastic) waves, in which case~$u,G,\rho$ and~$f$ take respectively the roles of transverse displacement, shear modulus, mass density, and body force. 

Recalling the plane wave expansion approach~\cite{NHA15, NSK12}, we seek the Bloch-wave solutions of~\eqref{PDE}, namely  
\begin{equation}\label{bloch1}
u(\bx) ~=~ \tilde{u}(\bx) \hh e^ {i\bk\cdot\bx}, \quad \tilde{u}: \hh\text{$Y$-periodic}
\end{equation}
where~$\tilde{u}$ depends implicitly on~$\bk\in\mathbb{R}^d$ and~$\omega\in\mathbb{R}$ -- which are hereon assumed to be fixed. If the source term is conveniently taken in the form of a plane wave, namely~$f(\bx) = \tilde{f} e^ {i\bk\cdot\bx}$ where $\tilde{f}$ is a constant, (\ref{PDE}) reduces to 
\begin{eqnarray} \label{CellPDE}
&&-\omega^2 \rho(\bx) \tilde{u} - \nabla_{\!\bk} \!\cdot\! \big(G(\bx)\nabla_{\!\bk}\tilde{u}\big) ~=~ \tilde{f} \quad \mbox{in~~} Y, \\
&& \begin{array}{rcl}
\tilde{u}|_{x_j=0} &\!\!\!=\!\!\!& \tilde{u}|_{x_j=\ell_j} \\*[2mm]
\bnu\cdot G\hh\nabla_{\!\bk}\tilde{u} |_{x_j=0} &\!\!\!=\!\!\!& -\bnu\cdot G\hh\nabla_{\!\bk}\tilde{u}|_{x_j=\ell_j}
\end{array}, \quad j=\overline{1,d}  \label{CellPDEBC1}
\end{eqnarray}
where $\nabla_{\!\bk}= \nabla \!+ i \bk$, $\,x_j=\bx\!\cdot\!\be_j$, and $\bnu$ is the unit outward normal on~$\partial Y$. 

\begin{remark} In the sequel, we implicitly assume that~$j=\overline{1,d}$ unless stated otherwise. \end{remark}

\subsection{Eigensystem representation} \label{Representation}

\noindent In what follows, we make use of the periodic function spaces
\begin{eqnarray*}
L^2_\rho(Y)  &\!\!=\!\!&  \{g\in L^2(Y):\, (g,\rho g) < \infty, \,\,g|_{x_j=0} = g|_{x_j=\ell_j}\}, \\
H^1_{p} (Y) &\! \!=\!\!& \{g\in L^2(Y):\, \nabla g \in (L^2(Y))^d, \,\,g|_{x_j=0} = g|_{x_j=\ell_j}\},
\end{eqnarray*}
where $(g,h)$ is the usual inner product for $g,h\in L^2(Y)$. To facilitate the treatment of vector- and tensor-valued functions, we shall also write 
\begin{equation}\label{innerprodY}
(\boldsymbol{g},\boldsymbol{h}) \;=\; \int_Y \boldsymbol{g}(\bx) : \bar{\boldsymbol{h}}(\bx)  \, \textrm{d}\bx, \qquad \boldsymbol{g} \in (L^2(Y))^{d^q}, \quad \boldsymbol{h} \in (L^2(Y))^{d^r}, \quad q,r\geqslant 0
\end{equation}
where ``:'' denotes the usual product, the inner product, and the $\min(q,r)$-tuple tensor contraction when $q=r=0$, $q=r=1$, and $\max(q,r)>1$, respectively. 

To cater for the asymptotic treatment of~\eqref{CellPDE}, let  $\tw \in H^1_p(Y)$ solve  
\begin{eqnarray} \label{CellPDEfw}
&& -\omega^2 \rho(\bx) \tw(\bx) - \nabla_{\!\bk}\!\cdot\!  \left( G(\bx) \nabla_{\!\bk} \tw \right) ~=~ 1 \quad \text{in}~~ Y, \\*[0.7mm]
&& \quad \bnu\cdot G \nabla_{\!\bk}\tw|_{x_j=0} ~=\; -\bnu\cdot G\nabla_{\!\bk}\tw|_{x_j=\ell_j}. \label{CellPDEfw2}
\end{eqnarray}
Hereon we denote by~$\bI$ the second-order identity tensor and, assuming Einstein summation notation, we let \mbox{$\nabla_{\!\bk}\hspace*{0.5pt} \bg \,=\, \be_j\otimes \partial\bg/\partial{x}_j +i\bk \otimes \bg\,$} for any vector or tensor field~$\bg$. 

\subsubsection{Eigensystem for the unit cell of periodicity}\label{SectionEigen}

\noindent From the variational formulation we find that $\left( - \nabla_{\!\bk} \!\cdot\! \left( G(\bx) \nabla_{\!\bk}  \right) \right)^{-1}$, as an operator from $L^2_\rho (Y)$ to itself with the range in $H^1_{p}(Y)$ subject to relevant boundary conditions, is a compact self-adjoint operator \cite{PBL78}. Hence for each~$\bk$ there exists an eigensystem $\{ \tilde{\phi}_{m} (\bk)\!\in\! H^1_{p}(Y), \,\tilde{\lambda}_{m}(\bk)\!\in\!\mathbb{R}\}$, conveniently normalized so that \mbox{$\|\tilde{\phi}_m\|_{L^2(Y)}=1$}, that solves
\begin{eqnarray} \label{Eigensystem}
&& -\nabla_{\!\bk} \!\cdot\! \big( G(\bx) \nabla_{\!\bk} \tilde{\phi}_{m} \big) ~=~ \tilde{\lambda}_{m}\rho \, \tilde{\phi}_m \quad \text{in}~~ Y, \\
&& ~~ \bnu\cdot G\nabla_{\!\bk}\tilde{\phi}_{m}|_{x_j=0} ~=\; -\bnu\cdot G\nabla_{\!\bk}\tilde{\phi}_m|_{x_j=\ell_j} \label{Eigensystem2}
\end{eqnarray}
where~$\{\tilde{\phi}_{m}\}$ are complete and orthogonal in $L^2_\rho(Y)$, namely  
\begin{eqnarray}\label{orthonormal}
(\rho\tilde{\phi}_m,\tilde{\phi}_n) \!\!\!&=&\!\!\! \delta_{mn} \hh (\rho\tilde{\phi}_n,\tilde{\phi}_n) \qquad \text{(no summation)}. 
\end{eqnarray}
With such definitions, the variational statement of~\eqref{CellPDE}--\eqref{CellPDEBC1} reads 
\begin{eqnarray} \label{CellPDEEigenstrainVariational}
-\omega^2 \int_Y \rho(\bx) \tilde{u}(\bx) \overline{ \tilde{\phi}}_m(\bx) \,\textrm{d}\bx \;-\; \int_Y \tilde{u}(\bx)\overline{\nabla_{\!\bk} \!\cdot\! \big( G(\bx) \nabla_{\!\bk}  \tilde{\phi}_m(\bx) \big) }\,\textrm{d}\bx ~=\;  
\int_Y \tilde{f} \, \overline{ \tilde{\phi}}_m(\bx) \,\textrm{d}\bx.  
\end{eqnarray}
By the completeness of $\{ \tilde{\phi}_m \}$ in $L^2_\rho (Y)$, every $\tilde{u}\in L^2_\rho (Y)$ can be written as
\begin{eqnarray} \label{eigen-exp1}
\tilde{u}(\bx) ~=~ \tilde{f}\hh\tw(\bx) ~=~ \tilde{f} \sum_{m=1}^\infty \tilde{\alpha}_m \tilde{\phi}_m(\bx), 
\end{eqnarray}
where $\tilde{\alpha}_m$ are constants. Thanks to~\eqref{CellPDEEigenstrainVariational} and the orthogonality of $\{\tilde{\phi}_m\}$ in $L^2_\rho(Y)$, we have 
\begin{eqnarray} \label{eigen-exp2}
\tilde{\alpha}_m ~=~ \frac{(1,\tilde{\phi}_m)}{(\rho\tilde{\phi}_m,\tilde{\phi}_m)} \hh
\frac{1}{\tilde{\lambda}_m \!-\omega^2} \quad \text{when} \quad \omega^2 \not= \tilde{\lambda}_m, \quad m\geqslant 1.
\end{eqnarray}


\subsection{Effective wave motion} \label{NonDegenerateCaseSubsection}

\noindent Let~$\mathcal{B}\ni\boldsymbol{0}$ denote the first Brillouin zone, i.e.~the counterpart of~$Y$ in the Fourier $\bk$-space containing the origin. With such definition, we concern ourselves with an effective description of the wave motion for~$\bk\in\mathcal{B}$. In the context of the Bloch wave equation~\eqref{CellPDE}, this is traditionally done by restricting~$(\bk,\omega)$ to some neighborhood of the first dispersion branch~$(\bk\in\mathcal{B},\omega=\tilde{\lambda}_1^{1/2}(\bk))$ and seeking the effective wave motion as $\langle\tilde{u}\rangle = (\tilde{u},1)$, i.e.~the $Y$-average of~$\tilde{u}$~\cite{NSK12} (recall that~$|Y|\!=\!1$). If further a \emph{low-wavenumber, low-frequency} (LW-LF) asymptotic model of the problem is sought, the effective i.e. homogenized field $\langle\tilde{u}\rangle$ is expanded about $(\bk=\bzero,\omega=0)$~\cite{MG17}.

Motivated by several recent developments~\cite{BS06,DBDF02,CKP10,HMC16}, we instead focus our efforts on formulating a comprehensive asymptotic framework toward \emph{low-wavenumber, finite-frequency} (LW-FF) approximation of~\eqref{CellPDE} that poses no upper limit on~$\omega$. When the effective behavior near the origin of~$\mathcal{B}$ is sought, this in particular suggests an expansion about \mbox{$(\bk=\bzero,\omega=\tilde{\lambda}_n^{1/2}(\bzero))$} for some $n> 1$. On recalling that in fact~$(\tilde{u},1)=(\tilde{u},\tilde{\phi}_{1}(\bzero))$, one may intuitively generalize upon the LW-LF definition of effective wave motion as
\begin{equation}\label{effective}
\langle\tilde{u}\rangle ~=~ (\tilde{u},\tilde{\phi}^{\circ}_{n}), \qquad \tilde{\phi}^{\circ}_{n}=\tilde{\phi}_n(\bzero), \qquad n\geqslant 1 
\end{equation}
to facilitate the homogenization at finite frequencies, see~\cite{CKP10} for example. In Section~\ref{main} we assume that the  eigenvalue~$\tilde{\lambda}_n(\bzero)$ is simple, and we use~\eqref{effective} as the basis for developing a second-order, LW-FF effective model of wave motion near~\mbox{$\bk=\bzero$}. In Section~\ref{Bri} and Section~\ref{Rep}, we extend the analysis to allow for \emph{finite-wavenumber, finite-frequency} (FW-FF) expansion about the ``corners''~$\bk^\ba$ of~$\mathcal{B}$, and situations when the eigenvalue~$\tilde{\lambda}_n(\bk^\ba)$ is repeated, respectively. 

With reference to~\eqref{effective}, we introduce the auxiliary ``zero-mean'' function space
\begin{equation}
H^1_{p0}(Y) ~=~ \{ g\in H^1_{p}(Y): \langle g\rangle=0\}, 
\end{equation}
where the choice of~$n\geqslant 1$ is implicit in the definition of~$\langle\cdot\rangle$. 

\begin{lemma} \label{Baiswv}
Let~$\tilde{u}$ solve~\eqref{CellPDE}--\eqref{CellPDEBC1} where~$\tilde{f}$ is a constant. Provided that $\omega^2 \neq \lambda_m \,\forall m$, we have 
\begin{eqnarray*}\label{effectivesol}
\tilde{u} ~=~ \tilde{f}\hspace*{0.3pt} \tw \qquad \text{and} \qquad 
\langle\tilde{u}\rangle ~=~ \tilde{f}\hspace*{0.3pt} \langle\tw\rangle.
\end{eqnarray*}
\end{lemma}

\section{LW-FF approximation near the origin of the first Brillouin zone} \label{main}

\noindent In principle, either the two-scale homogenization approach~\cite{WG15} or the asymptotically-expanded Willis' model~\cite{MG17} can be used to approximate the mean motion in a neighborhood of the \emph{acoustic branch}, $\omega=\tilde{\lambda}_1^{1/2}(\bk)$, at long wavelengths where $\|\bk\|\ll |Y|^{-1/d}$. With reference to the scalar wave equation~\eqref{PDE}, we seek the long-wavelength description of the effective field~$\langle u\rangle$ in a neighborhood of the $n$th branch, $\,\omega=\tilde{\lambda}_n^{1/2}(\bk)$, for $n\!\geqslant 1$. Accordingly, we describe the featured long-wavelength, finite-frequency (LW-FF) regime via scalings 
\begin{equation}\label{LW-FF}
\bk\,=\,\eps\hspace*{0.7pt} \hat{\bk}, \qquad \omega^2 \,=\, \tilde{\lambda}^{\circ}_n + \eps^2\nes\sigma\hat{\omega}^2, 
\qquad \tilde{\lambda}^{\circ}_n=\,\tilde{\lambda}_n(\bzero), \qquad \eps \,=\, o(1), \quad \sigma = \pm 1
\end{equation}
where~$\tilde{\lambda}^{\circ}_n\geqslant 0$ is such that 
\begin{eqnarray} \label{eignpde}
&& -\nabla \!\cdot\! \big(G\nabla \pho \big) ~=~ \tilde{\lambda}^{\circ}_n \hh\rho\hh \pho \quad \text{in}~~ Y, \\
&& \bnu\cdot G\nabla\pho|_{x_j=0} ~=\; -\bnu\cdot G\nabla\pho|_{x_j=\ell_j}, \label{eignbc}
\end{eqnarray}
and~$\pho\!\in\!H_p^1(Y)$ are orthogonal in $L^2_\rho(Y)$. For the time being, we assume that $\tilde{\lambda}^{\circ}_n$ to be \emph{simple}. 

\begin{remark}
Concerning the frequency scaling in~\eqref{LW-FF}, for~$n>1$ we have $\,\omega = \omega^{\circ}_n + \eps^2\nes\sigma \hat{\omega}^2/(2 \omega^{\circ}_n) +O(\eps^4)$, where $\omega^{\circ}_n=(\tilde{\lambda}^{\circ}_n )^{1/2}$. This is motivated by the observation that optical branches have zero initial slope, $\text{d}\omega/\text{d}\bk=\bzero$, at~$\bk=\bzero$ when~$\tilde{\lambda}^{\circ}_n$ is simple \cite{PBL78}, see also~\cite{CKP10} for the mathematical analysis in~$\mathbb{R}^2$. When~$\tilde{f}=0$ in~\eqref{CellPDE}, the sign factor~$\sigma$ in~\eqref{LW-FF} accordingly accounts for the possibility that the initial curvature of a given optical branch can be either positive or negative. When~$\tilde{f}\neq0$, on the other hand, the pair~$(\bk,\omega)$ is given a priori (say in a neighborhood of the~$n$th branch) and we take $\sigma=\text{sign}(\omega-\omega^{\circ}_n)$. 
\end{remark}

\begin{lemma} \label{lem2}
The solution of~\eqref{eignpde}--\eqref{eignbc} is characterized by~$\hh\arg(\pho(\bx))=\text{const.}\hh$ for~$\bx\in Y$. Accordingly we can take~$\pho(\bx)$ to be real-valued without loss of generality. See Appendix~A (electronic supplementary material) for proof.
\end{lemma}

\subsection{Asymptotic expansion at finite frequency} \label{aeei}

\noindent On recalling~\eqref{CellPDEfw}--\eqref{CellPDEfw2} and assuming the LW-FF regime~\eqref{LW-FF}, we find that $\tw\in H_p^1(Y)$ solves   
\begin{eqnarray} 
&&-(\tilde{\lambda}^{\circ}_n\!+\!\eps^2\nes\sigma\hat{\omega}^2) \rho\hh \tw -\big( \nabla \!+\! \eps\hh i \hat{\bk}\big) \!\cdot\! \big ( G ( \nabla \!+\! \eps\hh i \hat{\bk})  \tw \big) ~=~ 1 \quad \text{in}~~ Y, \label{ComparisonAnsatzEqn} \\ 
&& \qquad \bnu\cdot G(\nabla\!+\!\eps\hh i\hat{\bk})\tw|_{x_j=0} ~=\: -\bnu\cdot G(\nabla\!+\!\eps\hh i\hat{\bk})\tw|_{x_j=\ell_j}. \label{ComparisonAnsatzEqnBCN}
\end{eqnarray}
Consider next the asymptotic expansion 
\begin{eqnarray}\label{wexp}
\tw (\bx) ~=~ \eps^{-2} \tw_0(\bx) \;+\; \eps^{-1} \tw_1(\bx) \;+\; \tw_2(\bx) \;+\; \eps\hh \tw_3(\bx) \;+\; \cdots,
\end{eqnarray}
by which~(\ref{ComparisonAnsatzEqn})--(\ref{ComparisonAnsatzEqnBCN}) become a series in $\eps$. In what follows, the differential equations satisfied by $\tw_m$ in $Y$ ($m\geqslant 0$) are subject to implicit periodic boundary conditions
\begin{eqnarray}\label{ipbc}
\begin{array}{rcl} \tw_m|_{x_j=0} &\!\!\!=\!\!\!& \tw_m|_{x_j=\ell_j}, \\*[1.2mm]
\bnu\cdot G (\nabla\tw_{m}+i\hat{\bk}\hh\tw_{m-1})|_{x_j=0} &\!\!\!=\!\!\!& -\bnu\cdot G(\nabla\tw_{m}+ i\hat{\bk}\hh\tw_{m-1})|_{x_j=\ell_j}, 
\end{array} 
\end{eqnarray}
where $\tw_{-1}\equiv 0$. 

In the sequel, we conveniently denote by~$\ww_m$ the respective \emph{constants of integration} when solving for~$\tw_m(\bx)$. With reference to~\eqref{effective} and~\eqref{wexp}, we specifically seek~$\ww_m$ as    
\begin{equation}\label{effective2}
\ww_m ~=~ \langle\tw_m\rangle ~=~ (\tw_m,\pho), \qquad m\geqslant 0 
\end{equation}
so that 
\begin{equation}\label{wexpeff}
\langle\tw\rangle ~=~ \eps^{-2} \ww_0 \;+\; \eps^{-1} \ww_1 \;+\; \ww_2 \;+\; \eps\hh \ww_3 \;+\; \cdots.
\end{equation}

\subsubsection{Leading-order approximation}

\noindent The $O(\eps^{-2})$ contribution stemming from~\eqref{ComparisonAnsatzEqn} and~\eqref{wexp} reads 
\begin{eqnarray} \label{w0e}
&& -\tilde{\lambda}^{\circ}_n \hh\rho \hh\tw_0 -\nabla \!\cdot\! \big(G \nabla \tw_0 \big) ~=~ 0 \quad \text{in}~~ Y. 
\end{eqnarray}
By virtue of~\eqref{eignpde}--\eqref{eignbc} and~\eqref{ipbc} with~$m=0$, we have 
\begin{equation}\label{w0sol}
\tw_0(\bx) ~=~ \ww_0\hh \pho(\bx),
\end{equation} 
where $\ww_0$ is a constant to be determined. By the earlier premise that $\|\pho\|\!=\!1$, one immediately finds $\langle\tw_0\rangle = \ww_0$ as postulated by~\eqref{effective2}.

Similarly, the $O(\eps^{-1})$ equality can be identified as 
\begin{eqnarray}\label{w1e}
-\tilde{\lambda}^{\circ}_n\hh\rho\hh\tw_1
- \nabla \!\cdot\! \big(G (\nabla\tw_1 + i\hat{\bk}\hh\tw_0)\big) - i \hat{\bk}\hh \!\cdot\! \big( G \nabla \tw_0\big) ~=~ 0 \quad \text{in}~~ Y.
\end{eqnarray}
This equation is solved by
\begin{equation}\label{w1sol}
\tw_1(\bx) ~=~ \ww_0\hh \bchi^{\mbox{\tiny{(1)}}}(\bx)\!\cdot i\hat{\bk}  \;+\; \ww_1\hh\pho(\bx) 
\end{equation}
where~$\ww_1$ is a constant, and $\bchi^{\mbox{\tiny{(1)}}}\!\in(H_{p0}^1(Y))^d$ is a unique vector that satisfies 
\begin{eqnarray} \label{Comparisonchi1}
&& \tilde{\lambda}^{\circ}_n\hh\rho\hh\bchi^{\mbox{\tiny{(1)}}} + 
\nabla \!\cdot\! \big(G(\nabla\bchi^{\mbox{\tiny{(1)}}}\! + \bI\pho)\big) + G \nabla \pho ~=~0 \quad \text{in}~~ Y, \\
&& \bnu \cdot G(\nabla\bchi^{\mbox{\tiny{(1)}}}\! + \bI\pho)|_{x_j=0} ~=\: 
-\bnu \cdot G(\nabla\bchi^{\mbox{\tiny{(1)}}}\! + \bI\pho)|_{x_j=\ell_j}. \nonumber
\end{eqnarray}
Note that~\eqref{Comparisonchi1} is solvable since one can show that $\nabla \!\cdot\! (G \bI \pho)+ G \nabla \pho$ is orthogonal to $\pho$. In a similar fashion, the solvability condition can be demonstrated to hold for all subsequent cell problems and will not be discussed hereon. We also remark that the above expressions for $\tw_0$ and~$\tw_1$ are branch-generic, including the case of the acoustic branch ($n=1$) where $\tilde{\lambda}^{\circ}_1=0$ and~$\tilde{\phi}_1^\circ=1$, see~\cite{WG15, MG17} for example.

\begin{lemma} \label{lem2b}
Cell function~$\bchi^{\mbox{\tiny{(1)}}}$ is real-valued. See Appendix~A, electronic supplementary material, for proof.
\end{lemma} 

Proceeding with the cascade of differential equations, the $O(1)$ equation reads
\begin{eqnarray}\label{w2}
-\tilde{\lambda}^{\circ}_n\hh\rho\hh\tw_2
- \nabla \!\cdot\! \big(G (\nabla\tw_2 + i\hat{\bk}\hh\tw_1)\big) 
- i\hat{\bk} \!\cdot\! \big(G(\nabla\tw_1 + i \hat{\bk}\hh \tw_0)\big) 
- \sigma\hat{\omega}^2 \rho\hh \tw_0 ~=~ 1 \quad \text{in}~~ Y.  
\end{eqnarray}
To help expose the behavior of~$\ww_0$, we next evaluate the inner product of~\eqref{w2} with~$\pho$, i.e.
\begin{equation}\label{w2eff}
\big(\!-\tilde{\lambda}^{\circ}_n\hh\rho\hh\tw_2
- \nabla \!\cdot\! \big(G (\nabla\tw_2 \!+ i\hat{\bk}\hh\tw_1)\big), \pho\big) ~=~
\big(i\hat{\bk} \!\cdot\! \big(G(\nabla\tw_1 \!+ i \hat{\bk}\hh \tw_0)\big), \pho \big) + \sigma\hat{\omega}^2 \langle\rho\pho\rangle\hh \ww_0 + \langle 1\rangle.
\end{equation}
On deploying repeated integration by parts and recalling~\eqref{eignpde}, the second term on the left-hand side of~\eqref{w2eff} is computed as 
\begin{multline}
-\big(\nabla \!\cdot\! \big(G (\nabla\tw_2 + i\hat{\bk}\hh\tw_1)\big),\pho\big) ~=~ \big(G (\nabla\tw_2 + i\hat{\bk}\hh\tw_1)\hh ,\nabla \pho \big)  \\
~=~ - \big( \tw_2, \nabla \!\cdot\!  \big(G (\nabla\pho \big) \big) + \big(i\hat{\bk}G\tw_1 \hh ,\nabla \pho \big) 
~=~ \big(\tw_2, \tilde{\lambda}^{\circ}_n\hh\rho\hh\pho\big) + \big(i\hat{\bk}G\tw_1 \hh ,\nabla \pho \big). \label{lem3p2}
\end{multline}
By~\eqref{w1sol} and~\eqref{lem3p2}, the left-hand side of~\eqref{w2eff} becomes   
\[
\big(i\hat{\bk} G \tw_1,\nabla\pho\big) ~=~ \ww_0 \big(G\bchi^{\mbox{\tiny{(1)}}}\otimes\nabla\pho,1\big) : (i\hat{\bk})^2 + \ww_1\big(i\hat{\bk}\hh G\pho,\nabla\pho\big),
\]  
while its right-hand side reads 
\[
\ww_0  \big(G\nabla\bchi^{\mbox{\tiny{(1)}}},\pho\big) : (i\hat{\bk})^2 + \ww_1\big(i\hat{\bk} \!\cdot\!(G\nabla\pho),\pho\big) + \ww_0 \big(G\pho\bI,\pho\big) : (i\hat{\bk})^2 + \sigma\hat{\omega}^2 \langle\rho\pho\rangle\hh \ww_0 + \langle 1\rangle.
\]
From the last two results and the fact that~$\pho$ is real-valued, we find that~$\ww_0$ solves  
\begin{equation}\label{w0ode}
-\big(\bmu^{\mbox{\tiny{(0)}}}\!:(i\hat{\bk})^2 + \sigma\rho^{\mbox{\tiny{(0)}}}\hh\hat{\omega}^2\big)\hh \ww_0 ~=\; \langle1\rangle, 
\end{equation}
where the effective coefficients  
\begin{eqnarray} \label{effmod0}
\rho^{\mbox{\tiny{(0)}}} ~=~ \langle\rho\hh \pho\rangle, \qquad  
\bmu^{\mbox{\tiny{(0)}}} ~=~ \big\langle G\{\nabla\bchi^{\mbox{\tiny{(1)}}}+\bI\pho\}\big\rangle \,-\, \big (G\{\bchi^{\mbox{\tiny{(1)}}}\otimes\nabla\pho\},1 \big) 
\end{eqnarray}
are real-valued thanks to Lemma~\ref{lem2} and Lemma~\ref{lem2b}. In~\eqref{effmod0} and hereafter, $\{\boldsymbol{\cdot}\}$ denotes tensor averaging over all index permutations; in particular for an $n$th-order tensor~$\boldsymbol{\tau}$, one has 
\begin{equation}\notag
\{\boldsymbol{\tau}\}_{j_1,j_2,\ldots j_n} ~=~ \frac{1}{n!}\sum_{(l_1,l_2,\ldots l_n)\in P} \boldsymbol{\tau}_{l_1,l_2,\ldots l_n}, \qquad j_1,j_2,\ldots j_n \in\overline{1,d}
\end{equation}
where~$P$ denotes the set of all permutations of~$(j_1,j_2,\ldots j_n)$. Such averaged expression for~${\bmu}^{\mbox{\tiny{(0)}}}$ is due to the structure of~${\bmu}^{\mbox{\tiny{(0)}}} \!: (i\hat{\bk})^2$, which is invariant with respect to the index permutation of~${\bmu}^{\mbox{\tiny{(0)}}}$. Later on, we will also make use of the partial symmetrization 
\begin{equation}\notag
\{\boldsymbol{\tau}\}'_{j_1,j_2,\ldots j_n} ~=~ \frac{1}{(n\!-\!1)!}\sum_{(l_2,\ldots l_n)\in Q} \boldsymbol{\tau}_{j_1,l_2,\ldots l_n}, \qquad j_1,j_2,\ldots j_n \in\overline{1,d}
\end{equation}
where~$Q$ denotes the set of all permutations of~$(j_2,j_3,\ldots j_n)$. Note that assuming the acoustic branch ($\tilde{\phi}_1^\circ=1$) in~\eqref{effmod0}, one recovers the well-known result~$\rho^{\mbox{\tiny{(0)}}}=\langle\rho\rangle$ and $\bmu^{\mbox{\tiny{(0)}}}=\langle G\{\nabla\bchi^{\mbox{\tiny{(1)}}}+\bI\}\rangle$. 

\begin{remark}\label{rem2}
With reference to Lemma~\ref{Baiswv}, we first recall that~\eqref{w0ode} caters for an effective description of~\eqref{CellPDE} with~$\tilde{f}\neq 0$. In this case, we assume that ${\bmu}^{\mbox{\tiny{(0)}}} \!:\! (i\hat{\bk})^2  + \sigma \rho^{\mbox{\tiny{(0)}}}\hh\hat{\omega}^2 \neq 0$ so that~\eqref{w0ode} has a solution. By way of~\eqref{wexp} and~\eqref{w0sol}, the leading LW-FF approximation of the Bloch wave function~$\tilde{u}$ in~\eqref{bloch1} can thus be written as~$\tilde{u}_0 = \eps^{-2} \tilde{f}\tilde{w}_0 = \eps^{-2}\tilde{f}\ww_0\hh \pho$. Following the recent argument by Willis~\cite{W16}, we find that the mean (total) energy density of a Bloch wave~\eqref{bloch1} -- averaged in space over~$Y$ and time over~$2\pi/\omega$ -- is given by $\bar{E}=\tfrac{1}{2}\omega^2(\rho u, u) = \tfrac{1}{2}\omega^2 (\rho\tilde{u}, \tilde{u})$. Its leading LW-FF approximation thus reads $\bar{E}_0 = \tfrac{1}{2}\omega^2(\rho\tilde{u}_0, \tilde{u}_0) = \eps^{-4}\tfrac{1}{2}\rho^{\mbox{\tiny{(0)}}}\omega^2 |\tilde{f}\ww_0|^2$, which lends further credence to~\eqref{w0ode} as an effective descriptor of wave motion.  
\end{remark}

\begin{remark}\label{rem3}
Considering the free Bloch waves solving~\eqref{CellPDE} with~$\tilde{f}=0$, from~\eqref{w0ode} we directly obtain the leading LW-FF approximation of the $n$th dispersion branch, \mbox{$\omega(\bk)= \omega^{\circ}_n + \eps^2\nes\sigma \hat{\omega}^2(\hat{\bk})/(2 \omega^{\circ}_n)$}, by solving the primal problem 
\begin{equation}\label{disper0}
\bmu^{\mbox{\tiny{(0)}}}\!:(i\hat{\bk})^2 + \sigma\rho^{\mbox{\tiny{(0)}}}\hh\hat{\omega}^2=0. 
\end{equation}
To obtain a real-valued root for~$\hat{\omega}$, one must take $\sigma=sign(\bmu^{\mbox{\tiny{(0)}}})$, where $sign(\cdot)\in\{-1,1\}$ reflects the sign definiteness of its argument. 
\end{remark}

\subsubsection{First-order corrector}

\noindent Let $\bchi^{\mbox{\tiny{(2)}}} \in \big(H_{p0}^1(Y)\big)^{d \times d}$ be the unique second-order tensor solving 
\begin{eqnarray} \label{Comparisonchi2}
&& \hspace*{-7mm} \tilde{\lambda}^{\circ}_n\hh\rho\hh \bchi^{\mbox{\tiny{(2)}}} + 
\nabla \!\cdot\! \big( G\big( \nabla \bchi^{\mbox{\tiny{(2)}}} + \{\bI \otimes \bchi^{\mbox{\tiny{(1)}}}\}'\big)\big) + 
G\{\nabla\bchi^{\mbox{\tiny{(1)}}} + \bI\pho\} - \frac{\rho}{\rho^{\mbox{\tiny{(0)}}}}\hh \bmu^{\mbox{\tiny{(0)}}} \pho ~=~0 \quad \text{in}~~ Y, \\
&&  \bnu \cdot G\big( \nabla \bchi^{\mbox{\tiny{(2)}}} + \{\bI \otimes \bchi^{\mbox{\tiny{(1)}}}\}'  \big)|_{x_j=0} ~=\; - \bnu \cdot G\big( \nabla \bchi^{\mbox{\tiny{(2)}}} + \{\bI \otimes \bchi^{\mbox{\tiny{(1)}}}\}' \big)|_{x_j=\ell_j},  \nonumber
\end{eqnarray}
and let $\eta^{\mbox{\tiny{(0)}}}\in H_{p0}^1(Y)$ be the unique function satisfying  
\begin{eqnarray} \label{Comparisoneta0}
&& \tilde{\lambda}^{\circ}_n\hh\rho\hh\eta^{\mbox{\tiny{(0)}}} + 
\nabla \!\cdot\! \big(G\nabla \eta^{\mbox{\tiny{(0)}}}\big) - \frac{\rho}{\rho^{\mbox{\tiny{(0)}}}}\langle 1\rangle\hh \pho + 1 ~=~ 0 
\quad \text{in}~~ Y, \\
&& \quad \bnu \cdot G \nabla  \eta^{\mbox{\tiny{(0)}}}|_{x_j=0} ~=\: - \bnu \cdot G \nabla  \eta^{\mbox{\tiny{(0)}}}|_{x_j=\ell_j}. \nonumber
\end{eqnarray}
With such definitions, one can show that~\eqref{w2} is solved by 
\begin{eqnarray}\label{w2sol}
\tw_2 (\bx) \,=\, \ww_0\hh \bchi^{\mbox{\tiny{(2)}}}(\bx):(i\hat{\bk})^{2} \,+\, 
\ww_1\hh \bchi^{\mbox{\tiny{(1)}}} (\bx) \cdot i\hat{\bk}  \,+\, w_2\hh\pho(\bx) \,+\, \eta^{\mbox{\tiny{(0)}}}(\bx).
\end{eqnarray}

\begin{lemma} \label{lem4}
The following identity holds
\begin{eqnarray*}
\langle G \nabla \eta^{\mbox{\tiny{(0)}}} \rangle - (G \eta^{\mbox{\tiny{(0)}}},\nabla\pho)  ~=~
\frac{\langle 1\rangle}{\rho^{\mbox{\tiny{(0)}}}} \brho^{(1)} - (\bchi^{\mbox{\tiny{(1)}}},1). 
\end{eqnarray*}
See Appendix~A, electronic supplementary material, for proof.
\end{lemma}

Proceeding with the asymptotic analysis, the $O(\eps)$ statement stemming from~\eqref{ComparisonAnsatzEqn} is found as 
\begin{eqnarray}\label{w3}
-\tilde{\lambda}^{\circ}_n\hh\rho\hh\tw_3
- \nabla \!\cdot\! \big(G (\nabla\tw_3 + i\hat{\bk}\hh\tw_2)\big) 
- i\hat{\bk} \!\cdot\! \big(G(\nabla\tw_2 + i \hat{\bk}\hh \tw_1)\big) 
- \sigma\hat{\omega}^2 \rho\hh \tw_1 ~=~ 0 \quad \text{in}~~ Y.  
\end{eqnarray}
To help compute~$\ww_1$, we next evaluate the inner product of~\eqref{w3} with~$\pho$, namely 
\begin{equation}\label{w3eff}
\big(\!-\tilde{\lambda}^{\circ}_n\hh\rho\hh\tw_3
- \nabla \!\cdot\! \big(G (\nabla\tw_3 \!+ i\hat{\bk}\hh\tw_2)\big), \pho\big) ~=~
\big(i\hat{\bk} \!\cdot\! \big(G(\nabla\tw_2 \!+ i \hat{\bk}\hh \tw_1)\big), \pho \big) + \sigma\hat{\omega}^2 \langle\rho\tw_1\rangle.
\end{equation}
Integrating~\eqref{w3eff} by parts twice and making use of~\eqref{w0sol}, \eqref{w1sol}, \eqref{w2sol} and Lemma~\ref{lem4}, we obtain the governing equation for~$w_1$ as 
\begin{eqnarray} \label{w1ode0}
-\big({\bmu}^{\mbox{\tiny{(0)}}}\!: (i \hat{\bk})^{2}+ \sigma \rho^{\mbox{\tiny{(0)}}}\hh\hat{\omega}^2 \big) w_1  \,-\,  
\big(\bmu^{\mbox{\tiny{(1)}}}\!: (i \hat{\bk})^{3}  + \sigma \brho^{\mbox{\tiny{(1)}}} \!\!\cdot\! i\hat{\bk} ~\hat{\omega}^2 \big)w_0 ~=~
\Big(\frac{\langle 1\rangle}{\rho^{\mbox{\tiny{(0)}}}} \brho^{(1)} - (\bchi^{\mbox{\tiny{(1)}}},1)\Big) \!\cdot i\hat{\bk}
\end{eqnarray}
where 
\begin{eqnarray} \label{effmod1}
\brho^{\mbox{\tiny{(1)}}}  = \langle \rho\hh \bchi^{\mbox{\tiny{(1)}}}  \rangle, \qquad   
\bmu^{\mbox{\tiny{(1)}}} \;=\; \big\langle G\{\nabla\bchi^{\mbox{\tiny{(2)}}} + \bI \!\otimes\! \bchi^{\mbox{\tiny{(1)}}}\}  \big\rangle - 
\big(G\{\bchi^{\mbox{\tiny{(2)}}}\otimes\nabla\pho\},1\big).
\end{eqnarray}

\begin{remark}
Thanks to the fact that $\bchi^{\mbox{\tiny{(1)}}}$ is real-valued, one can show by following the proof of Lemma~\ref{lem2b} that $\bchi^{\mbox{\tiny{(2)}}}$ is real-valued as well. As a result, $\brho^{\mbox{\tiny{(1)}}}\in\mathbb{R}^d$ and $\bmu^{\mbox{\tiny{(1)}}}\in\mathbb{R}^{d\times d\times d}$.
\end{remark}

\begin{lemma} \label{lem5}
The homogenization coefficients in~\eqref{effmod1} satisfy the identity 
\begin{eqnarray} \label{mu1-dentity}
\rho^{\mbox{\tiny{(0)}}} \bmu^{\mbox{\tiny{(1)}}} ~=\;  \{\brho^{\mbox{\tiny{(1)}}}\! \otimes \bmu^{\mbox{\tiny{(0)}}}\}, 
\end{eqnarray} 
which reduces~\eqref{w1ode0} to 
\begin{eqnarray} \label{w1ode}
-\big(\bmu^{\mbox{\tiny{(0)}}} \!: (i \hat{\bk})^{2}+ \sigma \rho^{\mbox{\tiny{(0)}}}\hh\hat{\omega}^2 \big) w_1 ~=~ -(\bchi^{\mbox{\tiny{(1)}}},1) \cdot i\hat{\bk}.
\end{eqnarray}
See Appendix~A (electronic supplementary material) for proof.
\end{lemma}

\subsubsection{Second-order corrector}

\noindent Let $\bchi^{\mbox{\tiny{(3)}}} \in \big( H_{p0}^1(Y)\big)^{d\times d \times d}$ be the unique ``zero-mean'' third-order tensor solving 
\begin{eqnarray} \label{Comparisonchi3}
\hspace*{-8mm} && 
\tilde{\lambda}^{\circ}_n\hh\rho\hh \bchi^{\mbox{\tiny{(3)}}} + 
\nabla \!\cdot\! \big(G\big(\nabla \bchi^{\mbox{\tiny{(3)}}} + \{\bI \otimes \bchi^{\mbox{\tiny{(2)}}}\}' \big)\big)  + 
G \{\nabla \bchi^{\mbox{\tiny{(2)}}} \!+ \bI \otimes \bchi^{\mbox{\tiny{(1)}}}\} -\frac{\rho}{\rho^{\mbox{\tiny{(0)}}}}\{\bmu^{\mbox{\tiny{(0)}}} \!\otimes \bchi^{\mbox{\tiny{(1)}}}\} ~=~ 0 \quad \text{in} ~~ Y, \qquad \\
\hspace*{-8mm}&& \qquad\quad  \bnu \cdot G\big( \nabla \bchi^{\mbox{\tiny{(3)}}} + \{\bI \otimes \bchi^{\mbox{\tiny{(2)}}}\}'  \big)|_{x_j=0} ~=\: - \bnu \cdot G\big( \nabla \bchi^{\mbox{\tiny{(3)}}} + \{\bI \otimes \bchi^{\mbox{\tiny{(2)}}}\}'  \big)|_{x_j=\ell_j}, \nonumber
\end{eqnarray}
and let $\boldeta^{\mbox{\tiny{(1)}}} \in \big( H_{p0}^1(Y)\big)^{d}$ be the unique ``zero-mean'' vector given by 
\begin{eqnarray} \label{Comparisoneta1}
\hspace*{-8mm} &&\tilde{\lambda}^{\circ}_n\hh\rho\hh\boldeta^{\mbox{\tiny{(1)}}} + \nabla \!\cdot\! \big( G (\nabla \boldeta^{\mbox{\tiny{(1)}}} + \bI \eta^{\mbox{\tiny{(0)}}} ) \big) + G \nabla \eta^{\mbox{\tiny{(0)}}} + 
\frac{\rho}{\rho^{\mbox{\tiny{(0)}}}} \big((\bchi^{\mbox{\tiny{(1)}}},1)\pho - \langle 1\rangle \bchi^{\mbox{\tiny{(1)}}} \big) ~=~ 0 \quad \text{in}~~ Y, \\
\hspace*{-8mm} && \qquad \bnu \cdot G(\nabla \boldeta^{\mbox{\tiny{(1)}}} + \bI \eta^{\mbox{\tiny{(0)}}} )|_{x_j=0} ~=\: - \bnu \cdot G(\nabla \boldeta^{\mbox{\tiny{(1)}}} + \bI \eta^{\mbox{\tiny{(0)}}} )|_{x_j=\ell_j}. \nonumber
\end{eqnarray}
By virtue of~\eqref{Comparisonchi3} and~\eqref{Comparisoneta1}, \eqref{w3} is solved by 
\begin{equation}\label{w3sol}
\tw_3 (\bx) ~=~ \ww_0\hh \bchi^{\mbox{\tiny{(3)}}}(\bx):(i\hat{\bk})^3 \,+\, \ww_1\hh \bchi^{\mbox{\tiny{(2)}}}(\bx):(i\hat{\bk})^2  
\,+\, \ww_2\hh \bchi^{\mbox{\tiny{(1)}}}(\bx)\cdot i\hat{\bk} \,+\, \ww_3\hh\pho(\bx) \,+\, \boldeta^{\mbox{\tiny{(1)}}} (\bx) \cdot i \hat{\bk}. 
\end{equation} 

To compute~$\ww_2$, we next recall the $O(\eps^2)$ contribution to~\eqref{ComparisonAnsatzEqn}, i.e. 
\begin{eqnarray}\label{w4}
-\tilde{\lambda}^{\circ}_n\hh\rho\hh\tw_4
- \nabla \!\cdot\! \big(G (\nabla\tw_4 + i\hat{\bk}\hh\tw_3)\big) 
- i\hat{\bk} \!\cdot\! \big(G(\nabla\tw_3 + i \hat{\bk}\hh \tw_2)\big) 
- \sigma\hat{\omega}^2 \rho\hh \tw_2 ~=~ 0 \quad \text{in}~~ Y.  
\end{eqnarray}
Averaging this result in the sense of~\eqref{effective} yields the equation for constant~$w_2$ as 
\begin{multline} \label{w2ode}
-\big({\bmu}^{\mbox{\tiny{(0)}}} \!: (i\hat{\bk})^2 + \sigma\rho^{\mbox{\tiny{(0)}}}\hh\hat{\omega}^2 \big) w_2 - 
\big({\bmu}^{\mbox{\tiny{(2)}}} \!: (i\hat{\bk})^4 + \sigma\brho^{\mbox{\tiny{(2)}}} \!: (i\hat{\bk})^2\hat{\omega}^2\big) w_0  ~\,=\,~  \sigma \langle\rho\hh\eta^{\mbox{\tiny{(0)}}}\rangle\hh \hat{\omega}^2 \\
+ \Big(\big\langle G\{\nabla \boldeta^{\mbox{\tiny{(1)}}} \!+ \bI\eta^{\mbox{\tiny{(0)}}}\} \big\rangle - \big(G\{\boldeta^{\mbox{\tiny{(1)}}}\!\otimes\nabla\pho\},1\big)
+\frac{1}{\rho^{\mbox{\tiny{(0)}}}}\{\brho^{\mbox{\tiny{(1)}}}\!\otimes(\bchi^{\mbox{\tiny{(1)}}},1)\}  
\Big): (i\hat{\bk})^2,
\end{multline}
where 
\begin{eqnarray} \label{effmod2}
\brho^{\mbox{\tiny{(2)}}}  = \langle \rho\hh \bchi^{\mbox{\tiny{(2)}}}  \rangle, \qquad   
\bmu^{\mbox{\tiny{(2)}}} \;=\; \big\langle G\{\nabla\bchi^{\mbox{\tiny{(3)}}} + \bI \!\otimes\! \bchi^{\mbox{\tiny{(2)}}}\}  \big\rangle - 
\big(G\{\bchi^{\mbox{\tiny{(3)}}}\otimes\nabla\pho\},1\big).
\end{eqnarray}

\begin{remark}
In line with earlier arguments, we find that $\bchi^{\mbox{\tiny{(3)}}}\!:Y\to \mathbb{R}^{d\times d\times d}$, $\eta^{\mbox{\tiny{(0)}}}\!:Y\to\mathbb{R}$, $\boldeta^{\mbox{\tiny{(1)}}}\!:Y\to\mathbb{R}^{d}$, \mbox{$\brho^{\mbox{\tiny{(2)}}}\in\mathbb{R}^{d\times d}$}, and $\bmu^{\mbox{\tiny{(2)}}}\in\mathbb{R}^{d\times d\times d\times d}$.
\end{remark}

\begin{theorem} \label{explain1}
Consider the effective wave motion~$\langle\tilde{u}\rangle$ in the sense of~\eqref{effective}, where the Bloch wave function~$\tilde{u}$ solves~\eqref{CellPDE} with~$\tilde{f}\neq 0$. Assuming that the eigenvalue $\tilde{\lambda}^{\circ}_n$ ($n\geqslant1$) in~\eqref{eignpde} has multiplicity one, the second-order LW-FF approximation of~$\langle\tilde{u}\rangle$  in a neighborhood~\eqref{LW-FF} of~$(\bzero,\omega_n^\circ)$ satisfies
\begin{equation} \label{ode}
-\big({\bmu}^{\mbox{\tiny{(0)}}}\!: (i\hat{\bk})^{2}+ \sigma \rho^{\mbox{\tiny{(0)}}}\hh\hat{\omega}^2\big) \langle\tilde{u}\rangle \,-\,  
\eps^2\hh \big({\bmu}^{\mbox{\tiny{(2)}}} \!: (i\hat{\bk})^4 + \sigma\brho^{\mbox{\tiny{(2)}}} \!: (i\hat{\bk})^2\hat{\omega}^2\big) \langle\tilde{u}\rangle ~\,\overset{\eps}{=}\,~ \eps^{-2} \tilde{f}\hh M(\hat{\bk},\hat{\omega}),  
\end{equation}
where~``$\overset{\eps}{=}$'' signifies equality with an~$O(\eps)$ residual, and
\begin{multline} \label{M}
M(\hat{\bk},\hat{\omega}) ~=~ \langle 1 \rangle \,-\, \eps\, (\bchi^{\mbox{\tiny{(1)}}},1) \cdot i\hat{\bk} \,+\, \eps^2\hh \sigma \langle\rho\hh\eta^{\mbox{\tiny{(0)}}}\rangle\hh \hat{\omega}^2 \\*[-0.6mm]
+ \eps^2\hh \Big(\big\langle G\{\nabla \boldeta^{\mbox{\tiny{(1)}}} \!+ \bI\eta^{\mbox{\tiny{(0)}}}\} \big\rangle - \big(G\{\boldeta^{\mbox{\tiny{(1)}}}\!\otimes\nabla\pho\},1\big)
+\frac{1}{\rho^{\mbox{\tiny{(0)}}}}\{\brho^{\mbox{\tiny{(1)}}}\!\otimes(\bchi^{\mbox{\tiny{(1)}}},1)\}  
\Big): (i\hat{\bk})^2.
\end{multline}
When considering the propagation of free waves ($\tilde{f}=0$), on the other hand, the second-order LW-FF approximation of the $n$th dispersion branch reads
\begin{equation} \label{disper2ndb}
\sigma \hat{\omega}^2 ~\overset{\eps^3}{=}~ \frac
{{\bmu}^{\mbox{\tiny{(0)}}}\!: \hat{\bk}^{2} - \eps^2{\bmu}^{\mbox{\tiny{(2)}}} \!: \hat{\bk}^4} 
{\rho^{\mbox{\tiny{(0)}}} - \eps^2\hh \brho^{\mbox{\tiny{(2)}}} \!: \hat{\bk}^2}, \qquad n\geqslant 1.  
\end{equation}

\begin{proof}
From Lemma~\ref{Baiswv} and~\eqref{wexpeff}, we see
that~$\langle\tilde{u}\rangle = \tilde{f}(\eps^{-2}\ww_0+\eps^{-1}\hh\ww_1+\ww_2)+O(\eps)$. On computing the weighed sum $\tilde{f}(\eps^{-2}\eqref{w0ode}+\eps^{-1}\eqref{w1ode}+\eqref{w2ode})$, we recover~\eqref{ode} and thus~\eqref{disper2ndb}. 
\end{proof}
\end{theorem} 

\section{FW-FF approximation near ``corners'' of the first Brillouin zone} \label{Bri}

\noindent In what follows, we focus our attention on obtaining an effective description of wave motion inside the first ``quadrant'' of the first Brillouin zone, $\mathcal{B}=\{(\bk: |k_j|<\pi/\ell_j\}$, namely $\mathcal{B}^+ = \{\bk: 0< k_j < \pi/\ell_j\}$. The asymptotic treatment of $\mathcal{B}\backslash \overline{\mathcal{B}^+}$, as needed, can be performed in an analogous way. 

\subsection{Eigenfunction basis} \label{ebasis1}
 
\begin{figure}[h!] 
\centering{\includegraphics[width=\linewidth]{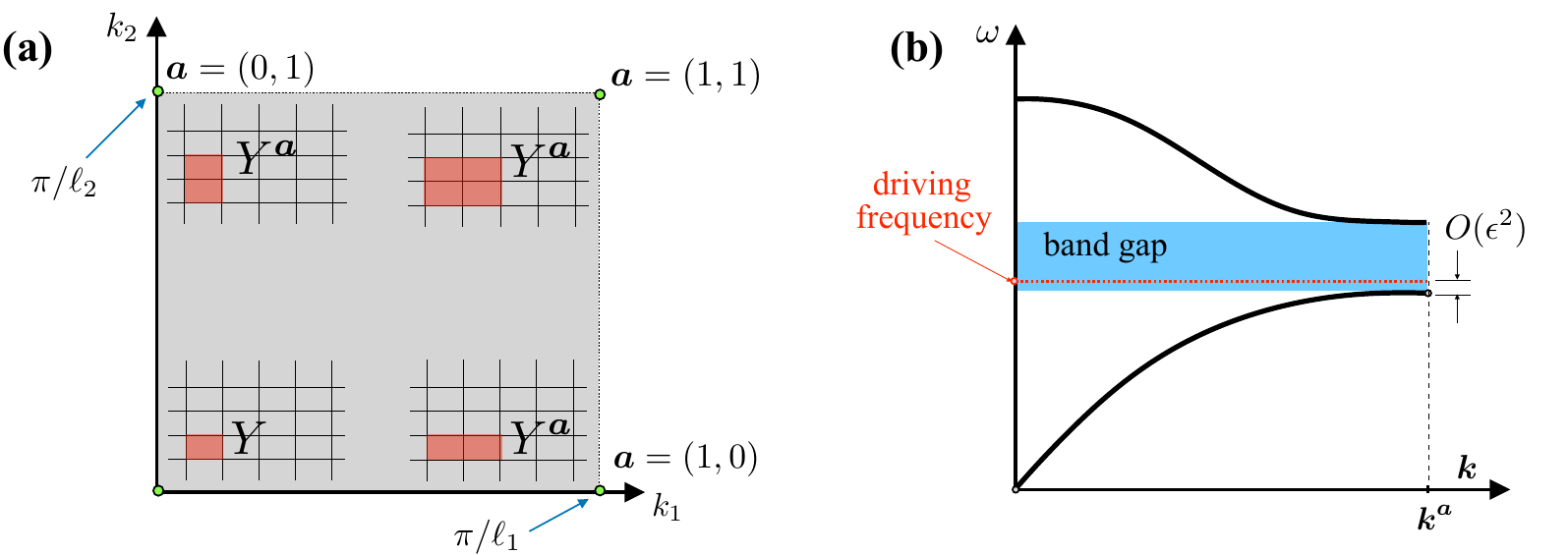}}\vspace*{-2mm}
\caption{Wave motion in a periodic medium: (a)~schematics of the multi-cell domains~$Y^\ba\!\in\mathbb{R}^2$ supporting the FW-FF homogenization of wave motion about vertices of the first quadrant of the first Brillouin zone, and (b) source excitation near the edge of a band gap.} \label{figo1} \vspace*{-2mm}
\end{figure}

\noindent With reference to~\eqref{Y}, consider the apexes of~$\mathcal{B}^+$ given by 
\begin{equation} \label{B+k}
\bk^{\ba} ~=~ \sum_{j=1}^{d} a_j \hh  \frac{\pi}{\ell_j} \be_j, \qquad \ba = (a_1,a_2,\ldots a_d), \qquad a_j\in\{0,1\}, 
\end{equation} 
where~$\be_j$ is the unit vector in the~$j$th coordinate direction.  At each such vertex, eigenfunction $\tilde{\phi}^\ba_{n}\in H^1_{p}(Y)$ corresponding to the $n$th solution branch satisfies 
\begin{eqnarray} \label{Eigensystem-apex11}
&&~~ -(\nabla\!+\!i\bk^\ba)\!\cdot\! \big(G(\nabla\!+\!i\bk^\ba)\tilde{\phi}^\ba_n\big) ~=~ \tilde{\lambda}_n^\ba \hh \rho \, \tilde{\phi}^\ba_n \quad \text{in}~~ Y, \\*[-0.5mm]
&& \bnu\cdot G(\nabla\!+\!i\bk^\ba)\tilde{\phi}^\ba_n|_{x_j=0} ~=\; -\bnu\cdot G(\nabla\!+\!i\bk^\ba)\tilde{\phi}^\ba_n|_{x_j=\ell_j}. 
\label{Eigensystem-apex12}
\end{eqnarray}
Letting 
\begin{equation}\label{Eigensystem-apex2}
\tilde{\phi}^\ba_n(\bx) ~=~ \tilde{\varphi}^\ba_n(\bx)\hh e^{-i\bk^\ba\!\cdot\bx} \quad \text{in}~~ Y, 
\end{equation} 
we find that 
\begin{eqnarray} \label{Eigensystem-apex3}
&& \hspace*{7mm} -\nabla\!\cdot\! \big(G\nabla\tilde{\varphi}^\ba_n\big) ~=~ \tilde{\lambda}_n^\ba \hh \rho \, \tilde{\varphi}^\ba_n \quad \text{in}~~ Y, \notag \\*[-0.8mm]
&& \begin{array}{rcl} 
\tilde{\varphi}^\ba_n|_{x_j=0} &\!\! \!=\!\!\!& (-1)^{a_j}\hh \tilde{\varphi}^\ba_n|_{x_j=\ell_j}, \\*[1.2mm]
\bnu\cdot G\nabla\tilde{\varphi}^\ba_n|_{x_j=0} &\!\! \!=\!\!\!& (-1)^{1+a_j}\hh\bnu\cdot G\nabla\tilde{\varphi}^\ba_n|_{x_j=\ell_j}, 
\end{array} 
\end{eqnarray}
which demonstrates that~$\tilde{\varphi}^\ba_n\,$ is $\,Y$-``anti-periodic'' (resp. $Y$-periodic) in the~$j$th coordinate direction when $a_j=1$ (resp.~$a_j=0$). Motivated by the multi-field continuum theory for two-dimensional crystal lattices~\cite{Vasi2006}, we next introduce the \emph{multi-cell domains}
\begin{equation}\label{Ya}
Y_\ba \,=\, \{\bx\!:\,0< \bx\!\cdot\!\be_j < (1\!+\!a_j)\ell_j\}, \qquad  |Y_\ba| \;=\: \prod_{j=1}^{d} (1+a_j), 
\end{equation}
and we extend the domain of~$\tilde{\varphi}^\ba_n\hh$ to $\hh Y_\ba$ as  
\begin{eqnarray} \notag
\forall \bx'\!\in Y, \quad \forall T\subseteq \{1,2,\ldots d\} \quad \Longrightarrow \quad 
\tilde{\varphi}^\ba_n (\bx) ~=~ (-1)^{\sum_{j\in T} a_j}\hh \tilde{\varphi}^\ba_n (\bx'), \quad  
\bx\,=\,\bx'+\sum_{j\in T} a_j \ell_j \be_j \in Y_\ba.
\end{eqnarray}
With such mapping, from~\eqref{Eigensystem-apex3} we find that~$\tilde{\varphi}^\ba_n \in H^1_{p}(Y_\ba)$ satisfy 
\begin{eqnarray} \label{Eigensystem-apex41}
&& \hspace*{4mm} -\nabla\!\cdot\! \big(G\nabla\tilde{\varphi}^\ba_n\big) ~=~ \tilde{\lambda}_n^\ba \hh \rho \, \tilde{\varphi}^\ba_n \quad \text{in}~~ Y_\ba, \\*[-0.8mm]
&& \begin{array}{rcl} 
a_j \tilde{\varphi}^\ba_n|_{x_j=0} &\!\!\!=\!\!\!& -a_j \tilde{\varphi}^\ba_n|_{x_j=\ell_j}, \\*[1.2mm] 
\bnu\cdot G\nabla\tilde{\varphi}^\ba_n|_{x_j=0} &\!\!\!=\!\!\!& -\bnu\cdot G\nabla\tilde{\varphi}^\ba_n|_{x_j=(1+a_j)\ell_j}, 
\end{array} \label{Eigensystem-apex43} 
\end{eqnarray}
where~$\bnu$ is the unit outward normal on~$\partial Y_\ba$. For clarity, domains~$Y_\ba$ are illustrated in Fig.~\ref{figo1}(a) for $d=2$. We note by recalling~\eqref{Eigensystem-apex2} and the normalization of~$\tilde\phi_n$ that~$(\pha,\pha)=(\vpha,\vpha)=1$. By virtue of this result and~\eqref{orthonormal}, we find~$\vpha$ to satisfy 
\begin{eqnarray} 
(\vpha,\vpha)_{Y_\ba} \;=\; |Y_\ba|, \qquad 
(\rho\tilde{\varphi}_m,\tilde{\varphi}_n)_{Y_\ba}  \:=\; \delta_{mn} \hh (\rho\tilde{\varphi}_n,\tilde{\varphi}_n)_{Y_\ba}  \quad \textrm{(no summation)}. \label{vorthonormal}
\end{eqnarray}
where~$(\cdot,\cdot)_{Y_\ba}$ denotes the inner product~\eqref{innerprodY} over~$Y_\ba$.

\begin{remark} \label{remfw1}
Thanks to~\eqref{Eigensystem-apex41}--\eqref{Eigensystem-apex43} and Lemma~\ref{lem2} by which~$\tilde{\varphi}^\ba_n$ can be taken as real-valued, it is apparent that the eigenfunction~$\tilde{\phi}^\ba_n$ corresponding to apex~$\ba$ of~$\mathcal{B}^+$ represents a plane wave propagating the stencil~$\tilde{\varphi}^\ba_n$ in direction~$-\bk^\ba$. When~$\bk^\ba=\bzero$, $\, \pha=\vpha\equiv\pho$ is a standing wave. 
\end{remark}

We will consider the effective wave motion near vertex~$\bk^\ba$ in terms of the averaging operators 
\begin{eqnarray}\label{effectivea}
\langle\tilde{u}\rangle_{\ba} &\!\!\!=\!\!\!& |Y_\ba|^{-1} \big(\tilde{u},\tilde{\phi}^{\ba}_{n}\big) _{Y_\ba}
 ~:=~  \big(\tilde{u},\tilde{\phi}^{\ba}_{n}\big)_{\overline{Y}_{\!\!\ba}}, \qquad n\geqslant 1, \\ \label{effectiveb}
\langle\tilde{g}\rangle_{\ba}^{\varphi} &\!\!\!=\!\!\!& |Y_\ba|^{-1} \big(\tilde{g},\tilde{\varphi}^{\ba}_{n}\big)_{Y_\ba} ~:=~ 
\big(\tilde{g},\tilde{\varphi}^{\ba}_{n}\big)_{\overline{Y}_{\!\!\ba}}, \qquad n\geqslant 1.
\end{eqnarray}
where~$\tilde{\phi}^{\ba}_{n}$ is extended to~$Y_\ba\supset Y$ by the application of~$Y$-periodicity. In the context of~\eqref{Eigensystem-apex2}, one sees that~\eqref{effectiveb} de facto ``averages'' its argument (by stencil~$\tilde{\varphi}^\ba_n$) over: (i) a single wavelength in the $\bk^\ba$-direction, and (ii)~dimension(s) of~$Y$ in the direction(s) perpendicular to~$\bk^\ba$. Note that in the LW-FF regime ($\ba=\bzero$), both~\eqref{effectivea} and~\eqref{effectiveb} recover~\eqref{effective}. In the context of~\eqref{effectiveb}, we also introduce the auxiliary function spaces
\begin{equation} \label{hilberta}
\begin{array}{rcl} 
H^{1}_{p}(Y_\ba) &\!\!\!=\!\!\!& \{ g\in L^2(Y_\ba): \, \nabla g \in (L^2(Y_\ba))^d, \,\,g|_{x_j=0} = g|_{x_j=(1+a_j)\ell_j}\},\\*[1.2mm] 
H^{1\ba}_{p0}(Y_\ba) &\!\!\!=\!\!\!&  \{ g\in H^{1}_{p}(Y_\ba): \langle g\rangle_{\ba}^{\varphi}=0\}. 
\end{array}
\end{equation}

\subsection{Scaling}\label{scaling}
 
\noindent For brevity of notation, we let~$\omega_n^\ba = (\tilde{\lambda}_n^\ba)^{1/2}$. Given the fact that the intersection between any flat portion of~$\partial\mathcal{B}^+$ and the $n$th dispersion branch is likely -- but not necessarily -- affiliated with the onset of a band gap, it is prudent to generalize~\eqref{LW-FF} as  
\begin{equation}\label{FW-FF}
\bk\,=\, \bk^\ba \!+ \eps\hh\hat{\bk}, \qquad \omega^2 \,=\, 
\tilde{\lambda}_n^\ba + \eps\hh\sigma\breve{\omega}^2 + \eps^2\nes\sigma\hat{\omega}^2, \qquad \eps \,=\, o(1), \quad \sigma = \pm 1. 
\end{equation}
On the basis of~\eqref{CellPDEfw}--\eqref{CellPDEfw2} and~\eqref{FW-FF}, we find that~$\tw\in H_p^1(Y)$ satisfies
\begin{eqnarray} 
&&-(\tilde{\lambda}^{\ba}_n\!+\! \eps\hh\sigma\hh\breve{\omega}^2\!+\!\eps^2\nes\sigma\hat{\omega}^2) \rho\hh \tw - \big(\nabla \!+\! i\bk^\ba \!+\! \eps\hh i \hat{\bk}\big) \!\cdot\! \big(G (\nabla \!+\! i\bk^\ba \!+\! \eps\hh i \hat{\bk}) \tw \big) ~=~ 1 \quad \text{in}~~ Y, \qquad  \label{ComparisonAnsatzEqna} \\*[-0.5mm] 
&&\quad \begin{array}{rcl} \tw|_{x_j=0} &\!\!=\!\!& \tw|_{x_j=\ell_j}, \\*[1.2mm]
\bnu\cdot G(\nabla \!+\! i\bk^\ba \!+\! \eps\hh i \hat{\bk})\tw|_{x_j=0} &\!\!=\!\!& -\bnu\cdot G(\nabla \!+\! i\bk^\ba \!+\! \eps\hh i \hat{\bk})\tw|_{x_j=\ell_j},
\end{array}\label{ComparisonAnsatzEqnBCNa}
\end{eqnarray}
where the first of~\eqref{ComparisonAnsatzEqnBCNa} is explicitly stated for reasons of convenience. 

To aid the asymptotic analysis, we next introduce~$\tw\in H_p^1(Y_a)$ as a field satisfying~\eqref{ComparisonAnsatzEqna}-\eqref{ComparisonAnsatzEqnBCNa} over~$Y_\ba$ by the application of $Y$-periodicity, and we consider the \emph{factorized} expansion 
\begin{eqnarray}\label{wexpa}
\tw (\bx) ~=~ e^{-i\bk^\ba\!\cdot\bx} \big(\eps^{-2} \tw_0(\bx) \;+\; \eps^{-1} \tw_1(\bx) \;+\; \tw_2(\bx) \;+\; \eps\hh \tw_3(\bx) \;+\; \cdots\big), 
\end{eqnarray}
which distills~\eqref{ComparisonAnsatzEqna}--\eqref{ComparisonAnsatzEqnBCNa} (repeated over~$Y_\ba$) into a series in $\eps$. By virtue of~\eqref{ComparisonAnsatzEqnBCNa}--\eqref{wexpa} and the analysis similar to that in Section~\ref{Bri}\ref{ebasis1}, we find that for every~$m\geqslant 0$, $\tw_m\!\in\!H_p^1(Y_\ba)$ are subject to the interior anti-symmetry conditions
\[
a_j \tw_m|_{x_j=0} ~=\; -a_j \tw_m|_{x_j=\ell_j},  
\] 
that are hereon assumed implicitly, and coupled boundary conditions 
\begin{eqnarray} \label{ipbca}
\bnu\cdot G (\nabla\tw_{m}+i\hat{\bk}\hh\tw_{m-1})|_{x_j=0} ~=\; -\bnu\cdot G(\nabla\tw_{m}+ i\hat{\bk}\hh\tw_{m-1})|_{x_j=(1+a_j)\ell_j},
\end{eqnarray}
where $\tw_{-1}\equiv 0$. In the sequel, we shall denote by~$\ww_m$ the respective weighted averages of~$\tw_m(\bx)$. Recalling~\eqref{effectivea}, \eqref{effectiveb} and~\eqref{wexpa}, we specifically seek~$\ww_m$ as    
\begin{equation}\label{effectivea2}
\ww_m ~=~ \langle\tw_m\rangle_{\ba}^{\varphi}, \qquad m\geqslant 0,
\end{equation}
so that~$\langle\tilde{u}\rangle_{\ba} = \tilde{f} \hh \langle\tw\rangle_{\ba}$ (by extension of Lemma~\ref{Baiswv}) and 
\begin{equation}\label{wexpeffa}
\langle\tw\rangle_{\ba} ~=~ \eps^{-2} \ww_0 \;+\; \eps^{-1} \ww_1 \;+\; \ww_2 \;+\; \eps\hh \ww_3 \;+\; \cdots.
\end{equation}

\subsection{Leading-order approximation} \label{loaa}

\noindent By the analyses given in Section~\ref{main}\ref{aeei} and Section~\ref{Bri}\ref{ebasis1}, the $O(\eps^{-2})$ contribution due to~\eqref{ComparisonAnsatzEqna} and~\eqref{wexpa} is found as  
\begin{eqnarray} \label{w0ea}
&& -\tilde{\lambda}^{\ba}_n \hh\rho \hh\tw_0 -\nabla\!\cdot\! \big(G \nabla \tw_0 \big) ~=~ 0 \quad \text{in}~~ Y_{\ba}. 
\end{eqnarray}
Due to~\eqref{Eigensystem-apex41}--\eqref{Eigensystem-apex43} and~\eqref{ipbca} with~$m=0$, we have 
\begin{equation}\label{w0sola}
\tw_0(\bx) ~=~ \ww_0\hh \tilde{\varphi}_n^\ba(\bx),
\end{equation} 
where $\ww_0$ is a constant to be determined. 

Similarly, the $O(\eps^{-1})$ equality can be identified as 
\begin{eqnarray}\label{w1ea}
-\tilde{\lambda}^{\ba}_n\hh\rho\hh\tw_1 - \sigma\breve{\omega}^2\!\rho\hh\tw_0
- \nabla \!\cdot\! \big(G(\nabla\tw_1 + i\hat{\bk}\hh\tw_0)\big) - i \hat{\bk}\hh \!\cdot\! \big( G \nabla \tw_0\big) ~=~ 0 \quad \text{in}~~ Y_{\ba}.
\end{eqnarray}
By the linearity of the problem, the solution of~\eqref{w1ea} is given by 
\begin{equation}\label{w1solaa}
\tw_1(\bx) ~=~ \ww_0 (\gamma^{\mbox{\tiny{(0)}}}(\bx)\hh \sigma\breve{\omega}^2 \!+ \bchi^{\mbox{\tiny{(1)}}}(\bx)\!\cdot i\hat{\bk})  \;+\; \ww_1\hh\vpha(\bx),
\end{equation}
where~$\ww_1$ is a constant; $\gamma^{\mbox{\tiny{(0)}}}\!\in H_{p0}^{1\ba}(Y_\ba)$ uniquely solves 
\begin{eqnarray} \label{Comparisoneta0a1}
&&~ \tilde{\lambda}^{\ba}_n\hh\rho\hh\gamma^{\mbox{\tiny{(0)}}} + 
\nabla \!\cdot\! \big(G\nabla \gamma^{\mbox{\tiny{(0)}}}\big) + \rho\hh\vpha ~=~ 0 \quad \text{in}~~ Y_{\ba}, \\
&& \bnu \cdot G \nabla \gamma^{\mbox{\tiny{(0)}}}|_{x_j=0} ~=\: - \bnu \cdot G \nabla  \gamma^{\mbox{\tiny{(0)}}}|_{x_j=(1+a_j)\ell_j}, \nonumber
\end{eqnarray}
and $\bchi^{\mbox{\tiny{(1)}}}\!\in(H_{p0}^{1\ba}(Y_\ba))^d$ is a unique vector that satisfies 
\begin{eqnarray} \label{Comparisonchi1a}
&& \quad \tilde{\lambda}^{\ba}_n\hh\rho\hh\bchi^{\mbox{\tiny{(1)}}} + 
\nabla \!\cdot\! \big(G(\nabla\bchi^{\mbox{\tiny{(1)}}}\! + \bI\vpha)\big) + G \nabla \vpha ~=~0 \quad \text{in}~~ Y_{\ba}, \\
&& \bnu \cdot G(\nabla\bchi^{\mbox{\tiny{(1)}}}\! + \bI\vpha)|_{x_j=0} ~=\: 
-\bnu \cdot G(\nabla\bchi^{\mbox{\tiny{(1)}}}\! + \bI\vpha)|_{x_j=(1+a_j)\ell_j}, \nonumber
\end{eqnarray}

\begin{lemma} \label{lem3}
Cell function~$\bchi^{\mbox{\tiny{(1)}}}$ solving~\eqref{Comparisonchi1a} is real-valued, and $\breve{\omega}=0$. See Appendix~A, electronic supplementary material, for proof.
\end{lemma} 

Thanks to Lemma~\ref{lem3}, the FW-FF asymptotic analysis near~$\bk=\bk^\ba$ largely reduces to its LW-FF companion near~$\bk=\bzero$ as in Section~\ref{main}. In particular, the $O(1)$ contribution to~\eqref{ComparisonAnsatzEqna} is found as 
\begin{equation}\label{w2a}
-\tilde{\lambda}^{\ba}_n\hh\rho\hh\tw_2 - \nabla \!\cdot\! \big(G (\nabla\tw_2 + i\hat{\bk}\hh\tw_1)\big) 
- i\hat{\bk} \!\cdot\! \big(G(\nabla\tw_1 + i \hat{\bk}\hh \tw_0)\big) - \sigma\hat{\omega}^2 \rho\hh \tw_0 ~=~ e^{i\bk^\ba\!\cdot\bx} \quad \text{in} ~~ Y_{\ba}.  
\end{equation}
By following the treatment of~\eqref{w2} that carries the same structure as~\eqref{w2a}, we find that~$\ww_0$ solves  
\begin{equation}\label{w0odea}
-\big(\bmu^{\mbox{\tiny{(0)}}}\!:(i\hat{\bk})^2 + \sigma\rho^{\mbox{\tiny{(0)}}}\hh\hat{\omega}^2\big)\hh \ww_0 ~=\; 
\langle e^{i\bk^\ba\!\cdot\bx}\rangle_{\ba}^{\varphi}, 
\end{equation}
featuring the real-valued effective coefficients  
\begin{eqnarray} \label{effmod0a}
\rho^{\mbox{\tiny{(0)}}} ~=~ \langle\rho\hh \vpha\rangle_{\ba}^{\varphi}, \qquad  
\bmu^{\mbox{\tiny{(0)}}} ~=~ \big\langle G\{\nabla\bchi^{\mbox{\tiny{(1)}}}+\bI\vpha\}\big\rangle_{\ba}^{\varphi} \,-\, 
\big (G\{\bchi^{\mbox{\tiny{(1)}}}\otimes\nabla\vpha\},1 \big)_{\overline{Y}_{\!\!\ba}}.
\end{eqnarray}

\begin{remark}
Assuming forced motion and taking $\bk$ and~$\omega$ as in~\eqref{FW-FF} with~$\breve{\omega}=0$ and~$n>1$, we see from Lemma~\ref{Baiswv}, \eqref{wexpa} and~\eqref{w0sola} that $u(\bx,t)=\tilde{u}(\bx)e^{i(\bk\cdot\bx-\omega t)}$ is computable as   
\begin{eqnarray} \notag
u(\bx,t) \!&=&\! \eps^{-2}\tilde{f} \hh \ww_0(\hat{\bk},\hat{\omega}) \, \vpha(\bx) e^{-i\bk^\ba\!\cdot\bx} \, 
e^{i((\bk^\ba+\eps\hat{\bk})\cdot\bx-[\omega_n^{\ba}+\eps^2\nes\sigma\hat{\omega}^2/(2\omega_n^{\ba})]t)} ~+~ O(\eps^{-1}), \quad \\ \notag
&=&\! \eps^{-2}\tilde{f} \hh \ww_0(\hat{\bk},\hat{\omega}) \hh \big[\vpha(\bx) \, e^{-i\omega_n^{\ba} t}\big] \, e^{i(\eps\hat{\bk}\cdot\bx -\eps^2\nes\sigma\hat{\omega}^2/(2\omega_n^{\ba})t)} ~+~ O(\eps^{-1}), 
\end{eqnarray}
where the bracketed term is a \emph{standing wave} over~$Y_\ba$. In the case of free Bloch waves, one similarly has 
\begin{eqnarray} \notag
u(\bx,t) \!&=&\! U \big[\vpha(\bx) \, e^{-i\omega_n^{\ba} t}\big] \,  e^{i(\eps\hat{\bk}\cdot\bx -\eps^2 \bmu^{\mbox{\tiny{(0)}}}:\hh\hat{\bk}^2/(2\rho^{\mbox{\tiny{(0)}}}\omega_n^{\ba})t)} + O(\eps), 
\end{eqnarray}
where~$U$ is a constant. This wave propagates the ``original'' Bloch eigenfunction~$\tilde{\phi}^\ba_n(\bx)$ with the phase and group velocities given respectively by 
\begin{eqnarray}  \notag
\boldsymbol{c}(\hat{\bk}) ~=~ 
\Big[\omega_n^{\ba}+ \eps^2  \frac{\bmu^{\mbox{\tiny{(0)}}}\!:\hat{\bk}^2}{2\rho^{\mbox{\tiny{(0)}}} \omega_n^{\ba}}\Big] \frac{\bk^\ba\!+\!\eps\hh\hat{\bk}}{\|\bk^\ba\!+\!\eps\hh\hat{\bk}\|^2}, \qquad\quad 
\boldsymbol{c}_g(\hat{\bk}) = \frac{\text{d}\omega}{\text{d}\bk}(\hat{\bk}) ~=~ \eps \, \frac{\bmu^{\mbox{\tiny{(0)}}}\!\cdot\hat{\bk}}{\rho^{\mbox{\tiny{(0)}}} \omega_n^{\ba}}.  
\end{eqnarray} 
\end{remark}

\subsection{First- and second-order correctors}\label{apex-2nd}

\noindent Proceeding with the analysis, we let $\bchi^{\mbox{\tiny{(2)}}} \in \big(H_{p0}^{1\ba}(Y_\ba)\big)^{d \times d}$ be the unique tensor solving 
\begin{eqnarray} \label{Comparisonchi2a}
&& \hspace*{-7mm} \tilde{\lambda}^{\ba}_n\hh\rho\hh \bchi^{\mbox{\tiny{(2)}}} + 
\nabla \!\cdot\! \big( G\big( \nabla \bchi^{\mbox{\tiny{(2)}}} + \{\bI \otimes \bchi^{\mbox{\tiny{(1)}}}\}'\big)\big) + 
G\{\nabla\bchi^{\mbox{\tiny{(1)}}} + \bI \vpha\} - \frac{\rho}{\rho^{\mbox{\tiny{(0)}}}} \bmu^{\mbox{\tiny{(0)}}} \vpha ~=~0 \quad \text{in}~~ Y_\ba, \\
&&  \bnu \cdot G\big( \nabla \bchi^{\mbox{\tiny{(2)}}} + \{\bI \otimes \bchi^{\mbox{\tiny{(1)}}}\}'  \big)|_{x_j=0} ~=\; - \bnu \cdot G\big( \nabla \bchi^{\mbox{\tiny{(2)}}} + \{\bI \otimes \bchi^{\mbox{\tiny{(1)}}}\}' \big)|_{x_j=(1+a_j)\ell_j},  \nonumber
\end{eqnarray}
and we introduce $\eta^{\mbox{\tiny{(0)}}}\in H_{p0}^{1\ba}(Y_\ba)$ as the unique function satisfying  
\begin{eqnarray} \label{Comparisoneta0a2}
&& \tilde{\lambda}^{\ba}_n\hh\rho\hh\eta^{\mbox{\tiny{(0)}}} + 
\nabla \!\cdot\! \big(G\nabla \eta^{\mbox{\tiny{(0)}}}\big) -
\frac{\rho}{\rho^{\mbox{\tiny{(0)}}}}\langle e^{i\bk^\ba\!\cdot\bx}\rangle_{\ba}^{\varphi} \, \vpha + e^{i\bk^\ba\!\cdot\bx} ~=~ 0 \quad \text{in}~~ Y_\ba, \\
&& \qquad \bnu \cdot G \nabla  \eta^{\mbox{\tiny{(0)}}}|_{x_j=0} ~=\: - \bnu \cdot G \nabla  \eta^{\mbox{\tiny{(0)}}}|_{x_j=(1+a_j)\ell_j}. \nonumber
\end{eqnarray}
Then one can show that 
\begin{eqnarray} \label{w1odea}
-\big(\bmu^{\mbox{\tiny{(0)}}} \!: (i \hat{\bk})^{2}+ \sigma \rho^{\mbox{\tiny{(0)}}}\hh\hat{\omega}^2 \big) w_1 ~=~ -(e^{i\bk^\ba\!\cdot\bx}\bchi^{\mbox{\tiny{(1)}}},1)_{\overline{Y}_{\!\!\ba}} \!\cdot i\hat{\bk}.
\end{eqnarray}

We next introduce $\bchi^{\mbox{\tiny{(3)}}} \in \big(H_{p0}^{1\ba}(Y_\ba)\big)^{d\times d \times d}$ as the unique ``zero-mean'' tensor solving 
\begin{eqnarray} \notag  
\hspace*{-5mm} && 
\tilde{\lambda}^{\ba}_n\hh\rho\hh \bchi^{\mbox{\tiny{(3)}}} + 
\nabla \!\cdot\! \big(G\big(\nabla \bchi^{\mbox{\tiny{(3)}}} \!+ \{\bI \otimes \bchi^{\mbox{\tiny{(2)}}}\}' \big)\big)  + 
G \{\nabla \bchi^{\mbox{\tiny{(2)}}} \!+ \bI \otimes \bchi^{\mbox{\tiny{(1)}}}\} -\frac{\rho}{\rho^{\mbox{\tiny{(0)}}}}\{\bmu^{\mbox{\tiny{(0)}}} \!\otimes \bchi^{\mbox{\tiny{(1)}}}\} ~=~ 0  
\quad \text{in} ~~ Y_\ba, \qquad \\
\hspace*{-5mm}&& \qquad  \bnu \cdot G\big( \nabla \bchi^{\mbox{\tiny{(3)}}} + \{\bI \otimes \bchi^{\mbox{\tiny{(2)}}}\}'  \big)|_{x_j=0} ~=\: - \bnu \cdot G\big( \nabla \bchi^{\mbox{\tiny{(3)}}} + \{\bI \otimes \bchi^{\mbox{\tiny{(2)}}}\}'  \big)|_{x_j=(1+a_j)\ell_j}, \nonumber
\end{eqnarray}
and we let $\boldeta^{\mbox{\tiny{(1)}}} \in \big(H_{p0}^{1\ba}(Y_\ba)\big)^{d}$ be the unique ``zero-mean'' vector given by 
\begin{eqnarray} \notag 
\hspace*{-7mm} && \tilde{\lambda}^{\ba}_n\hh\rho\hh\boldeta^{\mbox{\tiny{(1)}}} + \nabla \!\cdot\! \big( G (\nabla \boldeta^{\mbox{\tiny{(1)}}} \!+ \bI \eta^{\mbox{\tiny{(0)}}} ) \big) + G \nabla \eta^{\mbox{\tiny{(0)}}} + 
\frac{\rho}{\rho^{\mbox{\tiny{(0)}}}} \big((e^{i\bk^\ba\!\cdot\bx}\bchi^{\mbox{\tiny{(1)}}},1)_{\overline{Y}_{\!\!\ba}}\hh \vpha - 
\langle e^{i\bk^\ba\!\cdot\bx}\rangle_{\ba}^{\varphi}\hh \bchi^{\mbox{\tiny{(1)}}} \big) ~=~ 0  \quad \text{in}~~ Y_\ba, \qquad \\
\hspace*{-7mm} &&  \hspace*{10mm} \bnu \cdot G(\nabla \boldeta^{\mbox{\tiny{(1)}}} + \bI \eta^{\mbox{\tiny{(0)}}} )|_{x_j=0} ~=\: - \bnu \cdot G(\nabla \boldeta^{\mbox{\tiny{(1)}}} + \bI \eta^{\mbox{\tiny{(0)}}} )|_{x_j=(1+a_j)\ell_j}. \nonumber
\end{eqnarray}
With such definitions, one finds that 
\begin{multline} \label{w2odea}
-\big({\bmu}^{\mbox{\tiny{(0)}}} \!: (i\hat{\bk})^2 + \sigma\rho^{\mbox{\tiny{(0)}}}\hh\hat{\omega}^2 \big) w_2 - 
\big({\bmu}^{\mbox{\tiny{(2)}}} \!: (i\hat{\bk})^4 + \sigma\brho^{\mbox{\tiny{(2)}}} \!: (i\hat{\bk})^2\hat{\omega}^2\big) w_0  ~\,=\,~  
\sigma \langle\rho\hh\eta^{\mbox{\tiny{(0)}}}\rangle_{\ba}^\varphi\, \hat{\omega}^2 \\
+ \Big(\big\langle G\{\nabla \boldeta^{\mbox{\tiny{(1)}}} \!+ \bI\eta^{\mbox{\tiny{(0)}}}\} \big\rangle_{\ba}^\varphi - \big(G\{\boldeta^{\mbox{\tiny{(1)}}}\!\otimes\nabla\vpha\},1\big)_{\overline{Y}_{\!\!\ba}}
+\frac{1}{\rho^{\mbox{\tiny{(0)}}}}\big\{\brho^{\mbox{\tiny{(1)}}}\!\otimes (e^{i\bk^\ba\!\cdot\bx}\bchi^{\mbox{\tiny{(1)}}},1)_{\overline{Y}_{\!\!\ba}}\big\}  
\Big): (i\hat{\bk})^2,
\end{multline}
where 
\begin{eqnarray} \label{effmod2a}
\brho^{\mbox{\tiny{(2)}}}  = \langle \rho\hh \bchi^{\mbox{\tiny{(2)}}}  \rangle_{\ba}^{\varphi}, \qquad   
\bmu^{\mbox{\tiny{(2)}}} \;=\; \big\langle G\{\nabla\bchi^{\mbox{\tiny{(3)}}} + \bI \!\otimes\! \bchi^{\mbox{\tiny{(2)}}}\}  \big\rangle_{\ba}^{\varphi} - 
\big(G\{\bchi^{\mbox{\tiny{(3)}}}\otimes\nabla\vpha\},1\big)_{\overline{Y}_{\!\!\ba}}.
\end{eqnarray}

\begin{theorem} \label{explain1a}
Assume that the Bloch wave function~$\tilde{u}$ solves~\eqref{CellPDE} with~$\tilde{f}\neq 0$, and consider the effective wave motion $\langle\tilde{u}\rangle_{\ba}$ in the sense of~\eqref{effectivea} near apex~$\bk^\ba$ of the first ``quadrant'' of the first Brillouin zone~$\mathcal{B}^+$ according to~\eqref{B+k}. Assuming that the eigenvalue $\tilde{\lambda}^{\ba}_n$ ($n\!\geqslant\!1$) solving~\eqref{Eigensystem-apex41}--\eqref{Eigensystem-apex43} has multiplicity one, the second-order FW-FF approximation of~$\langle\tilde{u}\rangle_{\ba}$  in a neighborhood~\eqref{FW-FF} of~$\hh(\bk^\ba,\omega_n^\ba)$ satisfies
\begin{equation} \label{pde-fw-ff}
-\big({\bmu}^{\mbox{\tiny{(0)}}}\!: (i\hat{\bk})^{2}+ \sigma \rho^{\mbox{\tiny{(0)}}}\hh\hat{\omega}^2\big) \langle\tilde{u}\rangle_{\ba} \,-\,  
\eps^2\hh \big({\bmu}^{\mbox{\tiny{(2)}}} \!: (i\hat{\bk})^4 + \sigma\brho^{\mbox{\tiny{(2)}}} \!: (i\hat{\bk})^2\hat{\omega}^2\big) \langle\tilde{u}\rangle_{\ba} ~\,\overset{\eps}{=}\,~ \eps^{-2} \tilde{f}\hh M(\hat{\bk},\hat{\omega}),  
\end{equation}
where~``$\overset{\eps}{=}$'' signifies equality with an~$O(\eps)$ residual; the coefficients of homogenization $\rho^{\mbox{\tiny{(0)}}}\!\in\mathbb{R}, \hh {\bmu}^{\mbox{\tiny{(0)}}}\!\in\mathbb{R}^{d\times d}, \hh \brho^{\mbox{\tiny{(2)}}}\!\in\mathbb{R}^{d\times d}$ and~${\bmu}^{\mbox{\tiny{(2)}}}\!\in\mathbb{R}^{d\times d\times d\times d}$ are given by~\eqref{effmod0a} and~\eqref{effmod2a}; and 
\begin{multline}\notag 
M(\hat{\bk},\hat{\omega}) ~=~ \langle e^{i\bk^\ba\!\cdot\bx}\rangle_{\ba}^{\varphi} \,-\, 
\eps\, (e^{i\bk^\ba\!\cdot\bx}\bchi^{\mbox{\tiny{(1)}}},1)_{\overline{Y}_{\!\!\ba}} \!\cdot i\hat{\bk} \,+\, 
\eps^2\, \sigma \langle\rho\hh\eta^{\mbox{\tiny{(0)}}}\rangle_{\ba}^{\varphi} \, \hat{\omega}^2 \\ 
+ \eps^2\, \Big(\big\langle G\{\nabla \boldeta^{\mbox{\tiny{(1)}}} \!+ \bI\eta^{\mbox{\tiny{(0)}}}\} \big\rangle_{\ba}^{\varphi} - \big(G\{\boldeta^{\mbox{\tiny{(1)}}}\!\otimes\nabla\vpha\},1\big)_{\overline{Y}_{\!\!\ba}}
+\frac{1}{\rho^{\mbox{\tiny{(0)}}}}\big\{\brho^{\mbox{\tiny{(1)}}}\!\otimes (e^{i\bk^\ba\!\cdot\bx}\bchi^{\mbox{\tiny{(1)}}},1)_{\overline{Y}_{\!\!\ba}}\big\}  
\Big): (i\hat{\bk})^2.
\end{multline}
When considering the propagation of free waves ($\tilde{f}=0$), on the other hand, the second-order FW-FF approximation of the $n$th dispersion branch near~$\bk=\bk^\ba$ reads
\begin{equation}\notag
\sigma \hat{\omega}^2 ~\overset{\eps^3}{=}~ \frac
{{\bmu}^{\mbox{\tiny{(0)}}}\!: \hat{\bk}^{2} - \eps^2{\bmu}^{\mbox{\tiny{(2)}}} \!: \hat{\bk}^4} 
{\rho^{\mbox{\tiny{(0)}}} - \eps^2\hh \brho^{\mbox{\tiny{(2)}}} \!: \hat{\bk}^2}, \qquad n\geqslant 1.    
\end{equation}
\end{theorem}

\subsection{Example application: Green's function near the edge of a band gap}\label{space}

\noindent In the FW-FF regime covered by Theorem~\ref{explain1a}, the dominant wavelength ($2\pi/\|\bk^\ba\|$) is commensurate with the size of the unit cell, and the usual premise of the two-scale homogenization analysis ceases to apply. Instead, we make use of the Bloch expansion theorem~\cite{PBL78} to obtain an effective spatial description of wavefields whose spectrum is \emph{localized} in a neighborhood of~$(\bk^\ba,\omega_n^\ba)$. To facilitate the discussion, we temporarily adopt the terminology of the two-scale analysis and we let~$\by$ and~$\bx=\eps\by$ signify respectively the ``fast'' and ``slow'' spatial coordinate, where~$\by$ tracks fluctuations on the level of the unit cell. On denoting by~$\bgamma_\ba$ the lattice vectors of a structure with period~$Y_\ba$ and assuming that~$u$ solving~\eqref{PDE} for some~$f\!\in L^1(\mathbb{R}^d)$ is likewise absolutely integrable over~$\mathbb{R}^d$, their $Y_\ba$-periodic counterparts 
\begin{equation}\label{gfun0}
\tilde{g}(\by;\bk) \;= \sum_{\bgamma_\ba} g(\by+\bgamma_\ba) e^{-i\bk\cdot(\by+\bgamma_\ba)},  \qquad {g} = {u}, {f}
\end{equation}
can be shown to satisfy~\eqref{CellPDE}--\eqref{CellPDEBC1} with~$Y$ replaced by~$Y_\ba$. In this setting, the Bloch expansion theorem gives  
\begin{equation}\label{gfun1}
u(\by) \,=\, |\mathcal{B}_\ba|^{-1} \int_{\bk_s+\mathcal{B}_\ba} \tilde{u}(\by;\bk) \hh e^{i\bk\cdot\by}\hh \textrm{d}\bk, \qquad \by\in\mathbb{R}^d 
\end{equation}
and similarly for~$f$, where~$\mathcal{B}_\ba$ is the first Brillouin zone corresponding to~$Y_\ba$ and~$\bk_s\!\in\mathbb{R}^d$ is an arbitrary shift vector. In what follows, we conveniently take~$\bk_s=\bk^\ba$ and we denote by~$\bgamma$ the lattice vectors of a structure with period~$Y$.

We next consider the Green's function for an infinite periodic medium at frequency~$\omega$ near the edge of a band gap so that $\omega-\omega_n^\ba = O(\eps^2)$ for some~$n=\hat{n}$ and~$\bk=\bk^\ba$, whereas~$\omega-\omega_n(\bk)\geqslant O(1)$ for all other~$n$ and all~$\bk$ away from $\bk^\ba$ (see the schematics in Fig.~\ref{figo1}(b)). Letting the source density be~$f(\by)=\delta(\by-\by')$, one finds that the Green's function~$u$ under such conditions decays strongly~\cite{Doss2008} with~$\|\by-\by'\|$ so that~\eqref{gfun1} applies. 

By virtue of~\eqref{eigen-exp1}--\eqref{eigen-exp2}, the principal contribution to~\eqref{gfun1} in this case derives from a neighborhood of~$(\bk^\ba,\omega_{\hat{n}}^\ba)$. Assuming that~$\lambda_{\hat{n}}^\ba$ is an isolated eigenvalue, this neighborhood is characterized by scalings $\bk = \bk^\ba\!+ \eps\hat{\bk}$ and $\omega^2 = \lambda_{\hat{n}}^\ba + \eps^2\sigma\hat{\omega}^2$. From~\eqref{gfun0}, we find that~$\tilde{f}(\by;\bk) = \delta(\by\!-\!\by')e^{-i(\bk^\ba+\eps\hat\bk)\cdot\by}$ there. In this case ($\tilde{f}\neq \text{const.}$), the claim of Theorem~\ref{explain1a} can be generalized as 
\begin{eqnarray} \notag
-\big({\bmu}^{\mbox{\tiny{(0)}}}\!: (i\eps\hat{\bk})^{2}+ \rho^{\mbox{\tiny{(0)}}}\hh\eps^2 \sigma\hat{\omega}^2\big) \langle\tilde{u}\rangle_{\ba}  &\!\!\! = \!\!\!&
\langle e^{-i\eps\hat\bk\cdot\by}\delta(\by\!-\!\by')\rangle_{\ba}^{\varphi} \,-\, 
(e^{-i\eps\hat\bk\cdot\by}\bchi^{\mbox{\tiny{(1)}}}(\by) \delta(\by\!-\!\by'),1)_{\overline{Y}_{\!\!\ba}} \!\cdot i\eps\hat{\bk} \\ \label{gfun4}
&\!\!\! = \!\!\!& |Y_\ba|^{-1} \hh e^{-i\eps\hat\bk\cdot\by'} \big(\vphah(\by') \,-\, \bchi^{\mbox{\tiny{(1)}}}(\by')\cdot i\eps\hat{\bk}\big), 
\end{eqnarray}
upon discarding the second-order correction. 

By letting~$\tilde{\phi}^\ba_{\hat{n}}$ satisfy~\eqref{Eigensystem-apex11}--\eqref{Eigensystem-apex12} in every~$\bgamma+Y$, \eqref{gfun1} is next evaluated on the basis of~\eqref{gfun4} by: (i)~applying the local approximation 
\[
\tilde{u}(\by;\bk^\ba\!+\eps\hat{\bk}) ~\simeq~ \langle\tilde{u}\rangle_{\ba}[\eps\hat{\bk}] \, \hh\tilde{\phi}_{\hat{n}}^\ba(\by), \qquad \by\in\mathbb{R}^d 
\]
uniformly for all~$\eps\hat{\bk}\in\mathcal{B}_\ba$, and (ii) further extending the domain of integration to~$\mathbb{R}^d$, see~\cite{Doss2008} for related discussion. Specifically, we write 
\begin{equation}\label{gfun3}
u(\by) ~\overset{(i)}{\simeq}~ |\mathcal{B}_\ba|^{-1} \, \tilde{\phi}_{\hat{n}}^\ba(\by) \, e^{i\bk^\ba\cdot\by}  
\int_{\mathcal{B}_\ba} \!\! \langle\tilde{u}\rangle_{\ba}[\eps\hat{\bk}] \, e^{i\eps\hat\bk\cdot\by}\hh\hh \textrm{d}(\eps\hat\bk) ~=~ 
\tilde{\varphi}_{\hat{n}}^\ba(\by) \, U_\ba(\by), 
\end{equation}
where~$\tilde{\varphi}_{\hat{n}}^\ba=\tilde{\phi}_{\hat{n}}^\ba \, e^{i\bk^\ba\cdot\by}$ solves~\eqref{Eigensystem-apex41}--\eqref{Eigensystem-apex43} in every~$\bgamma_\ba\!+\!Y_\ba$, and 
\begin{equation}\label{gfun5}
U_\ba(\by) ~=~ |\mathcal{B}_\ba|^{-1} \int_{\mathcal{B}_\ba} \langle\tilde{u}\rangle_{\ba}[\eps\hat{\bk}] \, e^{i\eps\hat\bk\cdot\by}\hh\hh \textrm{d}(\eps\hat\bk)
~\overset{(ii)}{\simeq}~
|\mathcal{B}_\ba|^{-1} \int_{\mathbb{R}^d} \langle\tilde{u}\rangle_{\ba}[\eps\hat{\bk}] \, e^{i\eps\hat\bk\cdot\by}\hh\hh \textrm{d}(\eps\hat\bk).
\end{equation}
For brevity, we leave the error estimates for another study. We note, however, that the above approximations are supported by the facts that: (a)~$\langle\tilde{u}\rangle_{\ba}$ solving~\eqref{gfun4} behaves as $O((\eps^2\!+\|\eps\hat\bk\|^2)^{-1})$ for~$\eps\hat\bk\in\mathbb{R}^d$, and (b)~relative error due to approximation~(ii) diminishes fast with~$\|\by-\by'\|$ since the phase function in~\eqref{gfun5} has no stationary points. On recalling that~$|\mathcal{B}^\ba|=(2\pi)^d/|Y_\ba|$, \eqref{gfun4} and~\eqref{gfun5} demonstrate that~$U_\ba$ solves    
\begin{equation} \label{gfun6}
 -\big(\bmu^{\mbox{\tiny{(0)}}}\!\!:\!\nabla_{\!\by}^2 + \rho^{\mbox{\tiny{(0)}}}(\omega^2-\tilde{\lambda}_{\hat{n}}^\ba)\big) \hh U_{\ba} ~=~ 
\hh \big[\vphah(\by') \,-\, 
\bchi^{\mbox{\tiny{(1)}}}(\by') \cdot\! \nabla_{\!\by} \big]\hh F(\by-\by'),
\end{equation}
where 
\begin{equation} \label{gfun7}
F(\by-\by') ~=~ \left\{\begin{array}{ll}
|Y_\ba|^{-1} \prod_{j=1}^{d} \text{sinc}\big[\frac{\pi(y_j-y_j')}{(1+a_j)\ell_j}\big], \quad & \text{approx. (i)} \\*[1.5mm]
\delta(\by-\by'), & \text{approx. (ii)} 
\end{array}\right.  
\end{equation}
Later on, we will examine the utility of both approximations for a chessboard-like medium examined in~\cite{Cra11}. For future reference, we also note letting $d=2$, $\omega^2>\tilde{\lambda}_{\hat{n}}^\ba$, and~$\bmu^{\mbox{\tiny{(0)}}}=\mu^{\mbox{\tiny{(0)}}}\bI$ for some~$\mu^{\mbox{\tiny{(0)}}}\!<0$, that the solution of~\eqref{gfun7} affiliated with approximation~(ii) can be computed as 
\begin{equation} \label{gfun8}
 U_\ba^{\mbox{\tiny{(ii)}}}(\by) ~=~ \frac{1}{2\pi \mu^{\mbox{\tiny{(0)}}}} \Big[
\vphah(\by')\hh K_0(\alpha\|\by\|) \,+\, \frac{\alpha}{\|\by\|} \bchi^{\mbox{\tiny{(1)}}}(\by') \!\cdot\! \by \, K_1(\alpha\|\by\|)\Big],  
\qquad \alpha = \Big(\frac{ \rho^{\mbox{\tiny{(0)}}}(\omega^2-\tilde{\lambda}_{\hat{n}}^\ba)}{|\mu^{\mbox{\tiny{(0)}}}|}\Big)^{1/2},
\end{equation}
where~$K_0$ and~$K_1$ denote modified Bessel functions of the second kind. 

\section{Repeated eigenvalues}\label{Rep}

\noindent Up until this point, we assumed that the eigenvalue $\tilde{\lambda}^\ba_n$ is simple. In the case of repeated eigenvalues however, the asymptotic analysis must entail, and account for the interaction of, ``competing'' eigenfunctions at $\omega_n^\ba=(\tilde{\lambda}^\ba_n)^{1/2}$. This degenerate situation was recently examined in~\cite{CKP10,HMC16} and~\cite{Cra12} by considering one-dimensional wave motion, time-domain wave equation, and time-harmonic wave equation in separable bi-periodic structures, respectively. In composites with honeycomb structure, such situations are notably governed by the Dirac equation~\cite{TWZ17}. In this class of periodic configurations, salient features such as topologically protected edge states are possible -- leading to the concept of topological wave insulators, see~\cite{CKP10, HMC16, TWZ17} for details.

In what follows, we pursue the FW-FF expansion of wave motion about the apex point~$(\bk^\ba,\omega_n^\ba)$  in situations when the eigenvalue~$\lambda_n^\ba$ is repeated. For generality, we adopt the agnostic scaling framework~\eqref{FW-FF}, and we take~\eqref{ComparisonAnsatzEqna}-\eqref{ipbca} as the starting point.  

\subsection{Eigenfunction basis} \label{ebasis2}

\noindent Assuming that~$\lambda_n^\ba$ has multiplicity~$Q>1$, let $\tilde{\phi}^{\ba}_{nq}\in H_p^1(Y)$ $(q=1,Q)$ be linearly independent eigenfunctions that are orthogonal in $L^2_\rho(Y)$ and satisfy
\begin{eqnarray} \label{DoubleEigeignpde}
&&~~ -(\nabla\!+\!i\bk^\ba) \!\cdot\! \big(G (\nabla\!+\!i\bk^\ba) \tilde{\phi}^{\ba}_{nq} \big) ~=~ \tilde{\lambda}^{\ba}_n \hh\rho\hh \tilde{\phi}^{\ba}_{nq}  \quad \text{in}~~ Y, \\
&& \bnu\cdot G (\nabla\!+\!i\bk^\ba)\tilde{\phi}^{\ba}_{nq} |_{x_j=0} ~=\; -\bnu\cdot G (\nabla\!+\!i\bk^\ba)\tilde{\phi}^{\ba}_{nq} |_{x_j=\ell_j}. \label{Doubleeignbc}
\end{eqnarray}
Following the analysis in Section~\ref{Bri}\ref{ebasis1}, we let 
\begin{equation}\label{Eigensystem-apex22}
\tilde{\phi}^\ba_{nq}(\bx) ~=~ \tilde{\varphi}^\ba_{nq}(\bx)\hh e^{-i\bk^\ba\!\cdot\bx} \quad \text{in}~~ Y, \qquad q=\overline{1,Q} 
\end{equation} 
where~$\tilde{\varphi}^\ba_{nq} \in H^1_{p}(Y_\ba)$ are orthogonal according to~\eqref{vorthonormal} and solve
\begin{eqnarray} \label{Eigensystem-apex412}
&& \hspace*{5mm} -\nabla\!\cdot\! \big(G\nabla\tilde{\varphi}^\ba_{nq}\big) ~=~ \tilde{\lambda}_{n}^\ba \hh \rho \, \tilde{\varphi}^\ba_{nq} \quad \text{in}~~ Y_\ba, \\*[-0.8mm]
&& \begin{array}{rcl} \label{Eigensystem-apex432}
a_j \tilde{\varphi}^\ba_{nq}|_{x_j=0} &\!\!\!=\!\!\!& -a_j \tilde{\varphi}^\ba_{nq}|_{x_j=\ell_j},  \\*[1.2mm]
\bnu\cdot G\nabla\tilde{\varphi}^\ba_{nq}|_{x_j=0} &\!\!\!=\!\!\!& -\bnu\cdot G\nabla\tilde{\varphi}^\ba_{nq}|_{x_j=(1+a_j)\ell_j}.
\end{array}
\end{eqnarray}

By extending the argument of Lemma~\ref{lem2}, we find that~$\tilde{\varphi}^\ba_{nq}$ can be taken as real-valued without loss of generality. In this case, however, a generic solution to~\eqref{Eigensystem-apex412}--\eqref{Eigensystem-apex432} need not have constant argument across~$Y_\ba$.  

On the basis of~\eqref{DoubleEigeignpde}--\eqref{Eigensystem-apex432}, we define the relevant ``averaging'' operators as 
\begin{equation}\label{effectiveab2}
\langle\tilde{u}\rangle_{\ba}^{q} ~=~ 
\big(\tilde{u},\tilde{\phi}^{\ba}_{nq}\big)_{\overline{Y}_{\!\!\ba}}, \qquad
\langle\tilde{u}\rangle_{\ba}^{q\varphi} ~=~ 
\big(\tilde{u},\tilde{\varphi}^{\ba}_{nq}\big)_{\overline{Y}_{\!\!\ba}}, \qquad q=\overline{1,Q},
\end{equation}
and we generalize the second of~\eqref{hilberta} as  
\begin{equation} \label{hilbertap}
H^{1\ba}_{p0}(Y_\ba) ~=~  \{ g\in H^{1}_{p}(Y_\ba): \, \langle g\rangle_{\ba}^{q\varphi}=0, \, q=\overline{1,Q}\}.  
\end{equation}

\begin{remark} In the sequel, we do not asume implicit summation over repeated indexes~$p,q,r$ and~$s$. \end{remark}

\subsection{Leading-order approximation} \label{Loa}

\noindent As indicated earlier, we consider~\eqref{ComparisonAnsatzEqna}--\eqref{ComparisonAnsatzEqnBCNa} and we start from the asymptotic expansion~\eqref{wexpa} of~$\tw$. From the results in Section~\ref{Bri}\ref{loaa}, we find the $O(\eps^{-2})$ contribution stemming from~\eqref{ComparisonAnsatzEqna} and~\eqref{wexpa} to read as in~\eqref{w0ea}. By~\eqref{Eigensystem-apex412}--\eqref{Eigensystem-apex432}, we now have 
\begin{equation}\label{DoubleEigw0sol}
\tw_0(\bx) ~=~ \sum_q \ww_{0q} \hh \tilde{\varphi}^{\ba}_{nq}(\bx),
\end{equation} 
where $\ww_{0q}\!\in\mathbb{C}$ $(q=\overline{1,Q})$ are constants to be determined.

Next, we find the $O(\eps^{-1})$ statement of~\eqref{ComparisonAnsatzEqna} to be given by~\eqref{w1ea}. On multiplying the latter by~$\tilde{\varphi}^\ba_{np}$  and integrating over~$Y_\ba$, one can show by following the proof of Lemma~\ref{lem3} that 
\begin{eqnarray} \label{new1}
\sum_{q}\bth_{pq}^{\mbox{\tiny{(0)}}} \!\cdot\! (i\hat{\bk})\, \ww_{0q} \,+\,
\sigma \breve{\omega}^2 \rho_{p}^{\mbox{\tiny{(0)}}}\hh \ww_{0p} ~=~0, \qquad p=\overline{1,Q}. 
\end{eqnarray}
where 
\begin{equation} \label{pqmoduli}
\rho_{p}^{\mbox{\tiny{(0)}}} ~=~  \big\langle\rho\hh\tilde{\varphi}^{\ba}_{np}\big\rangle_{\ba}^{p\varphi}, \qquad 
\bth_{pq}^{\mbox{\tiny{(0)}}} ~=~ 
\big\langle G\nabla\tilde{\varphi}^{\ba}_{nq}\big\rangle_{\ba}^{p\varphi} - 
\big\langle G\nabla\tilde{\varphi}^{\ba}_{np}\big\rangle_{\ba}^{q\varphi}, \qquad p,q=\overline{1,Q}. 
\end{equation}

\begin{remark}
On recalling Lemma~\ref{Baiswv} and the discussion in Remark~\ref{rem2}, from~\eqref{DoubleEigw0sol} we find that the leading-order mean energy density of a Bloch wave~$\tilde{u}$ (averaged in space and time~\cite{W16}) is given by
\[
\bar{E}_0 ~=~ \tfrac{1}{2}\omega^2(\rho\tilde{u}_0, \tilde{u}_0) ~=~  \tfrac{1}{2}\omega^2 \tilde{f}^2 \sum_q \rho_{q}^{\mbox{\tiny{(0)}}} w_{0q}^2.  
\]
This sheds light on~$\rho_{q}^{\mbox{\tiny{(0)}}}$ ($q=\overline{1,Q}$) as the quanta of energy characterizing respective dispersion curves.
\end{remark}

For convenience, we rewrite~\eqref{new1} as the generalized eigenvalue problem
\begin{equation} \label{GEP}
\sum_{q} A_{pq} \hh\ww_{0q} \,-\, \tau \sum_{q} D_{pq} \hh\ww_{0q} ~=~0, \qquad p=\overline{1,Q} 
\end{equation}
with
\begin{equation} \label{ADmat}
A_{pq} ~=~ \bth_{pq}^{\mbox{\tiny{(0)}}} \!\cdot\! i\hat{\bk}, \qquad D_{pq} ~=~ \delta_{pq} \hh \rho_{q}^{\mbox{\tiny{(0)}}}, \qquad \tau  = -\sigma \breve{\omega}^2.  
\end{equation}
Since~$\bth_{pq}^{\mbox{\tiny{(0)}}}$ are real-valued and~$\rho_{q}^{\mbox{\tiny{(0)}}}>0$, we see that~$A_{pq}$ is Hermitian skew-symmetric, and that~$D_{pq}$ is positive definite. In what follows, we seek (whenever possible) a non-trivial solution to~\eqref{GEP}. 

\begin{lemma}\label{lem7}
The following statements hold: every~$\tau$ solving~\eqref{GEP} is either real-valued or zero; (ii) if~$\tau$ is an eigenvalue of~\eqref{GEP}, so is~$-\tau$; (iii) the eigenvectors corresponding to~$\tau$ and~$-\tau$ are complex conjugates of each other, and (iv) the maximum rank of~$A_{pq}$ is~$Q$ (resp.~$Q\!-\!1$) for~$Q$ even (resp.~odd). See Appendix~A (electronic supplementary material) for proof. 
\end{lemma}

\begin{assum} \label{fullrank}
Hereon, we focus our analysis on the situations where~$A_{pq}$ is either: (i)~of maximum rank according to Lemma~\ref{lem7}, or (ii) trivial.
\end{assum}

\begin{remark} 
The first step in the ensuing analyses is to expose the frequency scaling law in~\eqref{FW-FF} depending on the nature of~$A_{pq}$. Following the developments in Section~\ref{Bri} and in particular the proof of Lemma~\ref{lem3}, the key issue is whether letting~$\breve{\omega}=0$ turns the effective $O(\eps^{-1})$ equation~\eqref{GEP} into an identity for non-trivial~$\ww_{0q}$ -- which can then be deployed to satisfy the effective~$O(1)$ equation with the source term.  
\end{remark}


\subsubsection{Effective model for full-rank $A_{pq}$} \label{Apqnot0}

\noindent In situations when~$Q$ is even, generalized eigenvalue problem~\eqref{GEP} has \emph{no zero eigenvalues} by Lemma~\ref{lem7} and Assumption~\ref{fullrank}. Accordingly, assuming the scaling framework~\eqref{FW-FF} with~\mbox{$\breve{\omega}=0$} does not cater for a non-trivial solution of~\eqref{GEP}; we thus let~$\breve{\omega}\neq0$ and~$\ww_{0q}\!=0$ $(q=\overline{1,Q})$ by which the $O(\eps^{-1})$ contribution to~\eqref{ComparisonAnsatzEqna} becomes   
\begin{eqnarray}\label{even1}
-\tilde{\lambda}^{\ba}_n\hh\rho\hh\tw_1 - \nabla \!\cdot\! \big(G (\nabla\tw_1)\big) ~=~ 0 \quad \text{in}~~ Y_{\ba}.
\end{eqnarray}
As a result, the leading-order solution reads
\begin{equation}\label{even2}
\tw_1(\bx) ~=~ \sum_q \ww_{1q} \hh \tilde{\varphi}^{\ba}_{nq}(\bx),
\end{equation} 
where $\ww_{1q}\in\mathbb{C}$ are constants. To compute~$\ww_{1q}$, we identify the~$O(1)$ statement of~\eqref{ComparisonAnsatzEqna}, namely     
\begin{eqnarray}\label{even3}
-\tilde{\lambda}^{\ba}_n\hh\rho\hh\tw_2 - \nabla \!\cdot\! \big(G (\nabla\tw_2 
+ i\hat{\bk}\hh\tw_1)\big) - i \hat{\bk}\hh \!\cdot\! \big( G \nabla \tw_1\big) - \sigma \breve{\omega}^2\rho\hh\tw_1 ~=~ e^{i\bk^\ba\cdot\bx} \quad \text{in}~~ Y_\ba.
\end{eqnarray}
On integrating~$\langle\eqref{even3}\rangle_{\ba}^{p\varphi}$, $p\in\overline{1,Q}$ by parts as in the treatment of~\eqref{lem3p2}, we obtain
\begin{eqnarray} \label{even4}
\big(i\hat{\bk}\hh G \tw_1, \nabla\tilde{\varphi}^{\ba}_{np}\big)_{\overline{Y}_{\!\!\ba}} 
- \big(i\hat{\bk} \!\cdot\! \big( G \nabla \tw_1\big),\tilde{\varphi}^{\ba}_{np}\big)_{\overline{Y}_{\!\!\ba}}  
- \sigma \breve{\omega}^2 \big(\rho\hh\tw_1,\tilde{\varphi}^{\ba}_{np}\big)_{\overline{Y}_{\!\!\ba}} ~=~
\langle e^{i\bk^\ba\cdot\bx}\rangle_{\ba}^{p\varphi}. 
\end{eqnarray}
The substitution of~\eqref{even2} into~\eqref{even4} results in a system of~$Q$ governing equations for~$\ww_{1p}$ as 
\begin{eqnarray} \label{even5}
-\sum_{q} A_{pq}\hh \ww_{1q} \,-\, \sigma \breve{\omega}^2 \sum_{q} D_{pq}\hh \ww_{1q} ~=~ 
\langle e^{i\bk^\ba\cdot\bx}\rangle_{\ba}^{p\varphi}, \qquad p=\overline{1,Q} 
\end{eqnarray}
where~$A_{pq}$ and~$D_{pq}$ are given by~\eqref{ADmat}. By Lemma~\ref{lem7} and Assumption~\ref{fullrank}, we find that the generalized eigenvalue problem underpinning~\eqref{even5} has~$Q$ real eigenvalues, hereon denoted by~$-\sigma_r\hh\breve{\omega}_r^2$, and we assume~$\sigma\breve{\omega}^2 \neq \sigma_r\hh\breve{\omega}_r^2$ $(r=\overline{1,Q})$ when solving~\eqref{even5}.

\begin{remark} \label{rem12}
By Lemma~\ref{lem7}, the above eigenvalues appear in ``$\pm$'' pairs, and we conveniently arrange them so that (i)~$\breve{\omega}_{r}^2=\breve{\omega}_{r+1}^2$ and $\sigma_{r/r+1}\!=\!\pm 1$ for~$r$ odd, and (ii) $0< \breve{\omega}_1^2 \leqslant \breve{\omega}_2^2 \leqslant \ldots \leqslant \breve{\omega}_{Q}^2$. In this setting, there are~$Q/2$ pairs of dispersion branches emanating from the apex point~$(\bk^\ba,\omega_n^\ba)$ given by mappings \mbox{$(\bk^\ba\!+\eps\hh\hat{\bk}, \hh\omega_n^\ba\!+\eps\hh\sigma_r\breve{\omega}_r^2(\hat{\bk})/(2 \omega_n^\ba))$}. Each such branch is characterized by an~$O(1)$ (positive or negative) group velocity, and is stenciled by the respective eigenfunction~$\sum_q e^r_q\hh\vphaq$, where~$e^r_q\in\mathbb{C}$ ($q=\overline{1,Q}$) are components of the unit eigenvector of~$A_{pq}$ corresponding to~$\sigma_r\hh\breve{\omega}_r^2$. 
\end{remark}

\subsubsection{Effective model for trivial~$A_{pq}$}\label{apq0}

\noindent In certain situations, the ``primal'' matrix~$A_{pq}$ may become zero due to either trivial~$\bth_{pq}^{\mbox{\tiny{(0)}}}$ in~\eqref{pqmoduli}, or the vanishing products~$\bth_{pq}^{\mbox{\tiny{(0)}}}\!\cdot i\hat{\bk}$. In this case, (\ref{GEP}) with~$A_{pq}=0$ yields~$\breve{\omega}=0$ in order to preserve the leading-order solution~\eqref{DoubleEigw0sol}. Accordingly, from~\eqref{w1ea} with~$\breve{\omega}=0$ and~\eqref{even2} we find that  
\begin{equation}\label{w1solac}
\tw_1(\bx) ~=~ \sum_q \ww_{0q}\hh \bchi_q^{\mbox{\tiny{(1)}}}(\bx)\!\cdot i\hat{\bk}  \;+\; \sum_q \ww_{1q}\hh\vphaq(\bx),
\end{equation}
where~$\ww_{1q}\!\in\mathbb{C}$ $(q=\overline{1,Q})$ are constants, and $\bchi_q^{\mbox{\tiny{(1)}}}\in (H_{p0}^{1\ba}(Y_\ba))^d$ solve  
\begin{eqnarray} \label{Comparisonchi1ax}
\hspace*{-7mm} && \tilde{\lambda}^{\ba}_n\hh\rho\hh\bchi_q^{\mbox{\tiny{(1)}}} + 
\nabla \!\cdot\! \big(G(\nabla\bchi_q^{\mbox{\tiny{(1)}}}\! + \bI\vphaq)\big) + G \nabla \vphaq 
- \sum_r\frac{1}{\rho_{r}^{\mbox{\tiny{(0)}}}}\hh\bth_{rq}^{\mbox{\tiny{(0)}}}\hh\rho\hh\vphar ~=~0 \quad \text{in}~~ Y_{\ba}, \quad q=\overline{1,Q} \qquad \\*[-0.5mm] \notag
\hspace*{-7mm} && \quad\qquad \bnu \cdot G(\nabla\bchi_q^{\mbox{\tiny{(1)}}}\! + \bI\vphaq)|_{x_j=0} ~=\: 
-\bnu \cdot G(\nabla\bchi_q^{\mbox{\tiny{(1)}}}\! + \bI\vphaq)|_{x_j=(1+a_j)\ell_j}. 
\end{eqnarray}
Note that the last (summative) entry in~\eqref{Comparisonchi1ax} is inserted to ensure that the source term driving~$\bchi_q^{\mbox{\tiny{(1)}}}$ is orthogonal to~$\vphar$ $(r=\overline{1,Q})$, while making \emph{no contribution} to~\eqref{w1ea} due to the fact that~$A_{rq}=0$. On substituting this result into~\eqref{w2a}, multiplying by~$\tilde{\varphi}^\ba_{np}$, and integrating over~$Y_\ba$ we find that  
\begin{equation} \label{odd7A0}
-\sum_{q} B_{pq} \hh\ww_{0q} \,-\, \sigma\hat{\omega}^2  \sum_{q} D_{pq} \hh\ww_{0q} ~=~   
\langle e^{i\bk^\ba\cdot\bx}\rangle_{\ba}^{p\varphi}, \qquad p=\overline{1,Q},   
\end{equation}
where $D_{pq}$ is given by~\eqref{ADmat} and 
\begin{equation} \label{Bmat}
\hspace*{-7mm}
\qquad B_{pq} =\bmu_{pq}^{\mbox{\tiny{(0)}}}:(i\hat{\bk})^2, \qquad
\bmu_{pq}^{\mbox{\tiny{(0)}}} ~=~ \big\langle G\{\nabla\bchi_q^{\mbox{\tiny{(1)}}} +\bI \vphaq\}\big\rangle_{\ba}^{p\varphi} - 
\big(G\{\bchi_q^{\mbox{\tiny{(1)}}}\otimes\nabla\vphap\},1\big)_{\overline{Y}_{\!\ba}}.   
\end{equation}

\begin{lemma}\label{lem8}
Matrix $B_{pq}$ is symmetric. See Appendix~A, electronic supplementary material, for proof.
\end{lemma}

\subsubsection{Effective model for $rank(A_{pq})=Q\!-\!1$} \label{Apqneg1}

\noindent When~$Q$ is odd, by Lemma~\ref{lem7} and Assumption~\ref{fullrank} the generalized eigenvalue problem~\eqref{GEP} has \emph{one zero eigenvalue} to which corresponds unit eigenvector with components~${\vartheta}_q\in\mathbb{R}, \,q=\overline{1,Q}$. Without loss of generality, we select the eigenfunction basis~$\vphaq$ so that~$\vartheta_q=\delta_{q1}$, i.e. 
\begin{equation} \label{A000}
A_{1p}=A_{p1}=0, \quad p = \overline{1,Q}.
\end{equation}
In this case, \eqref{GEP} turns into identity if: (i)~$\hh\tau = -\sigma\breve{\omega}^2=0$, and (ii)~$\hh w_{0q} = \text{w}_0\hh \delta_{q1}$ i.e. 
\begin{equation} \label{new2}
\tw_{0}(\bx) ~=~ \text{w}_0 \hh \vphao(\bx), 
\end{equation}
for some constant~$\text{w}_0\in\mathbb{C}$. From~\eqref{w1ea} with~$\breve{\omega}=0$ and~\eqref{new2}, we find that  
\begin{equation}\label{w1solab}
\tw_1(\bx) ~=~ \text{w}_0 \bchi_{1}^{\mbox{\tiny{(1)}}}(\bx)\!\cdot i\hat{\bk}  \;+\; \sum_q \ww_{1q}\hh\vphaq(\bx)
\end{equation}
where~$\ww_{1q}\!\in\mathbb{C}$ ($q=\overline{1,Q}$) are constants, and $\,\bchi_{1}^{\mbox{\tiny{(1)}}}\!\in\!(H_{p0}^{1\ba}(Y_\ba))^d$ satisfies~\eqref{Comparisonchi1ax} with~$q=1$.

To compute~$\text{w}_0$, we recall the~$O(1)$ statement of~\eqref{ComparisonAnsatzEqna} given by~\eqref{w2a}, whose multiplication by~$\tilde{\varphi}^\ba_{n1}$ and integration over~$Y_\ba$ yields 
\begin{equation} \label{odd7}
- \big(\bmu_{1\nes 1}^{\mbox{\tiny{(0)}}}\!:\!(i\hat{\bk})^2 + 
\sigma\hat{\omega}^2\rho_{1}^{\mbox{\tiny{(0)}}}\big) \text{w}_0 ~=\, 
 \langle e^{i\bk^\ba\cdot\bx}\rangle_{\ba}^{1\varphi}, 
\end{equation}
where~$\rho_{1}^{\mbox{\tiny{(0)}}}$ and~$\bmu_{1\nes 1}^{\mbox{\tiny{(0)}}}$ are given respectively by~\eqref{pqmoduli} and~\eqref{Bmat}. On adopting the reasoning as in Section~\ref{main}, from~\eqref{odd7} one finds that~$\sigma=sign(\omega\!-\!\omega_n^\circ)$ when~$\tilde{f}\!\neq\!0$, and~$\sigma=sign(\bmu^{\mbox{\tiny{(0)}}})$ when~$\tilde{f}\!=\!0$. 

We remark however that the above leading-order model is \emph{incomplete} in that it features only a single projection, $\text{w}_0$, at an apex point that features~$Q$ eigenfunctions. Such description is indeed insufficient when~$ \langle e^{i\bk^\ba\cdot\bx}\rangle_{\ba}^{1\varphi}=0$. To identify the remaining~$Q-1$ leading-order projections in terms of~$w_{1q}$, we revisit the~$O(1)$ statement of~\eqref{ComparisonAnsatzEqna} taking~$\breve{\omega}\neq 0$, namely
\begin{equation} \notag 
-\tilde{\lambda}^{\ba}_n\hh\rho\hh\tw_2 - \nabla\!\cdot\!\big(G(\nabla\tw_2 + i\hat{\bk}\hh\tw_1)\big) - 
i \hat{\bk}\hh \!\cdot\! \big(G(\nabla\tw_1 + i\hat\bk\tw_0)\big) - \sigma\breve{\omega}^2\rho\hh\tw_1  - \sigma\hat{\omega}^2\rho\hh\tw_0 ~=~ e^{i\bk^\ba\cdot\bx} ~~ \text{in}~\; Y_\ba,
\end{equation}
where~$\tw_0$ is given by~\eqref{new2}. By following the analysis in Section~\ref{Rep}\ref{Loa}\ref{Apqnot0} and making use of~\eqref{Bmat}, \eqref{A000} and~\eqref{odd7}, we find that 
\begin{equation} \label{odd8}
-\rho_{1}^{\mbox{\tiny{(0)}}} w_{1\nes 1} ~=~  \text{w}_0  \langle \rho\hh\bchi_{1}^{\mbox{\tiny{(1)}}} \rangle_{\ba}^{1\varphi} \!\cdot\!i\hat{\bk},
\end{equation}
while~$w_{1q}$, $q=\overline{2,Q}$ solve the linear system 
\begin{multline} \label{odd9}
-\sum_{q>1} A_{pq}\hh \ww_{1q} \,-\, \sigma \breve{\omega}^2 \sum_{q>1} D_{pq}\hh \ww_{1q} ~=~ 
\langle e^{i\bk^\ba\cdot\bx}\rangle_{\ba}^{p\varphi} \:+\: 
\text{w}_0 \hh \big(\bmu_{p1}^{\mbox{\tiny{(0)}}}\!:\!(i\hat{\bk})^2 + \sigma\breve{\omega}^2 \brho_{p1}^{\mbox{\tiny{(1)}}} \!\cdot\!i\hat{\bk}\big), \\*[-2.7mm] 
\brho_{p1}^{\mbox{\tiny{(1)}}} = \langle\rho\hh\bchi_{1}^{\mbox{\tiny{(1)}}}\rangle_{\ba}^{p\varphi}, \qquad p=\overline{2,Q}.
\end{multline}
Note that the role of~$w_{1\nes 1}$ in the above expansion is ``secondary'' in that it vanishes when~$\text{w}_0\!=\!0$; a situation that motivates the consideration of~$w_{1q}$ in the first place. As a result the ``full'' effective model when~$rank(A_{pq})=Q\!-\!1$, that embodies the leading-order contribution of \emph{each eigenfunction}~$\vphaq$ ($q=\overline{1,Q}$), is thus given by~$\text{w}_0$ solving~\eqref{odd7} and~$w_{1q}$, $q=\overline{2,Q}$ solving~\eqref{odd9}. For completeness, we remark that the existence of this class of configurations (termed Dirac-like points) at~$\bk^\ba\!=\!\bzero$ with~$Q\!=\!3$ is directly related to the phenomenon of zero-index metamaterials~\cite{Ash2015}.  

\begin{remark}
On the basis of the above analysis, one may expect that the leading-order FW-FF approximation at some~$(\bk^\ba,\omega_n)$ with~$rank(A_{pq})=Q\!-\!N$, $N\!>\!1$ and~$Q\!-\!N$ even, will be described in terms of~$N$ components of~$w_{0q}$ and~$Q\!-\!N$ components of~$w_{1q}$. 
\end{remark} 

\subsubsection{Synthesis} \label{synthrepeig}

\begin{theorem} \label{61}
Consider the Bloch wave function~$\tilde{u}(\bx)$ solving~\eqref{CellPDE} with~$\tilde{f}\neq 0$ near apex~$\bk^\ba$ of the first ``quadrant'' of the first Brillouin zone, see~\eqref{B+k}. Assuming that the eigenvalue $\tilde{\lambda}^{\ba}_n>0$ solving~\eqref{Eigensystem-apex41}--\eqref{Eigensystem-apex43} has multiplicity~$Q>1$ and that the coupling matrix $A_{pq}$ in~\eqref{GEP} is either non-degenerate in the sense of Assumption~\ref{fullrank} or zero, the leading FW-FF approximation of $u(\bx,t)=\tilde{u}(\bx)\hh e^{i(\bk\cdot\bx-\omega t)}$ in a neighborhood~\eqref{FW-FF} of~$\hh(\bk^\ba,\omega_n^\ba)$ reads  
\begin{eqnarray*} 
u(\bx,t) &\!\!=\!\!& 
\left\{ \begin{array}{ll}
\eps^{-1}\tilde{f} \sum_r \ww_{1r}(\hat{\bk},\breve{\omega}) \hh \big[\vphar(\bx) \, e^{-i\omega_n^{\ba} t}\big] \, e^{i(\eps\hat{\bk}\cdot\bx -\eps\hh\sigma\breve{\omega}^2/(2\omega_n^{\ba})t)}, & rank(A_{pq})=Q, \\*[2mm]
\eps^{-2}\tilde{f} \sum_r \ww_{0r}(\hat{\bk},\hat{\omega}) \big[\vphar(\bx) \, e^{-i\omega_n^{\ba} t}\big] \, e^{i(\eps\hat{\bk}\cdot\bx -\eps^2\nes\sigma\hat{\omega}^2/(2\omega_n^{\ba})t)}, \quad & rank(A_{pq})=0, \\*[2mm]
\eps^{-2}\tilde{f} \text{\emph{w}}_0 (\hat{\bk},\breve{\omega}) \big[\vphao(\bx) \, e^{-i\omega_n^{\ba} t}\big] \, e^{i(\eps\hat{\bk}\cdot\bx -\eps^2\nes\sigma\hat{\omega}^2/(2\omega_n^{\ba})t)} \quad & \\*[0.4mm] 
\: + \: \eps^{-1}\tilde{f} \sum_s \ww_{1s}(\hat{\bk},\breve{\omega}) \hh \big[\vphas(\bx) \, e^{-i\omega_n^{\ba} t}\big] \, e^{i(\eps\hat{\bk}\cdot\bx -\eps\hh\sigma\breve{\omega}^2/(2\omega_n^{\ba})t)}, & rank(A_{pq})=Q\!-\!1, 
\end{array}\right. \quad 
\end{eqnarray*}
where~$r=\overline{1,Q}$; $s=\overline{2,Q}$; $\hh\ww_{1r}$ solve system~\eqref{even5}; $\ww_{0r}$ solve system~\eqref{odd7A0}; $\text{\emph{w}}_0$ solves~\eqref{odd7}, $\hh\ww_{1s}$ solve system~\eqref{odd9} and, when $rank(A_{pq})=Q\!-\!1$, the eigenfunction basis is chosen so that~\eqref{A000} holds. When considering the free waves~$(\tilde{f}=0)$, we have 
\begin{eqnarray*} 
u(\bx,t) &\!\!=\!\!& 
\left\{ \begin{array}{ll}
\sum_s U_s\hh \big[\sum_r e^s_r\hh \tilde{\varphi}_{nr}^\ba (\bx) \, e^{-i\omega_n^{\ba} t}\big] \, e^{i(\eps\hat{\bk}\cdot\bx -\eps\hh\sigma_s\breve{\omega}_s^2(\hat{\bk})/(2\omega_n^{\ba})t)}, &  rank(A_{pq})=Q,  \\*[2mm]
\sum_s U_s\hh \big[\sum_r g^s_r\hh \tilde{\varphi}_{nr}^\ba (\bx) \, e^{-i\omega_n^{\ba} t}\big] \, e^{i(\eps\hat{\bk}\cdot\bx -\eps^2\nes\sigma_s \hat{\omega}_s^2(\hat{\bk})/(2\omega_n^{\ba})t)}, & rank(A_{pq})=0,    \\*[2mm]
U_0\hh \big[\tilde{\varphi}_{n1}^\ba(\bx) \, e^{-i\omega_n^{\ba} t}\big] \, e^{i(\eps\hat{\bk}\cdot\bx -\eps^2\nes\sigma_0 \hat{\omega}_0^2(\hat{\bk})/(2\omega_n^{\ba})t)} & \\*[0.4mm]
\: + \: \sum_m U_m\hh \big[\sum_r h^m_r\hh \tilde{\varphi}_{nr}^\ba (\bx) \, e^{-i\omega_n^{\ba} t}\big] \, e^{i(\eps\hat{\bk}\cdot\bx -\eps\hh\sigma_m\breve{\omega}_m^2(\hat{\bk})/(2\omega_n^{\ba})t)}, &  rank(A_{pq})=Q\!-\!1, 
\end{array}\right.  
\end{eqnarray*}
where~$r,s\!=\!\overline{1,Q}$; $m\!=\!\overline{2,Q}$; $\hh U_0, U_s$ and~$U_m$ are constants; $-\sigma_s\hh\breve{\omega}_s^2(\hat{\bk})$ (resp.~$e^s_r$) are the eigenvalues (resp.~eigenvectors) of system~\eqref{even5}; $-\sigma_s\hh\hat{\omega}_s^2(\hat{\bk})$ (resp.~$g^s_r$) are the real eigenvalues (resp.~real eigenvectors) of system~\eqref{odd7A0}; $-\sigma_0 \hat{\omega}_0^2(\hat{\bk})$ is the eigenvalue of~\eqref{odd7}; and  $-\sigma_m\hh\breve{\omega}_m^2(\hat{\bk})$ (resp.~$h^m_r$) are the eigenvalues (resp.~eigenvectors) of system~\eqref{odd9}. The situations with $rank(A_{pq})=Q$ (resp.~$rank(A_{pq})=Q-1$) require that~$Q$ be even (resp.~odd). 
\end{theorem}

\begin{remark} \label{dirac}
From~\eqref{ADmat}, \eqref{even5}, \eqref{odd7A0}, and~\eqref{Bmat} we see that~$\breve{\omega}_s^2(\hat{\bk})\propto\|\hat{\bk}\|$ and~$\hat{\omega}_s^2(\hat{\bk})\propto\|\hat{\bk}\|^2$ in any given direction~$\hat\bk/\|\hat{\bk}\|$. Due to this result and~\eqref{FW-FF}, one finds that
\begin{equation}
\omega_s(\bk)\!-\!\omega_n^\ba ~\propto~ \left\{
\begin{array}{lll} 
\|\bk\!-\!\bk^\ba\| ~&\text{for}&~ rank(A_{pq})=Q, \\
\|\bk\!-\!\bk^\ba\|^2 ~&\text{for}&~ rank(A_{pq})=0, \end{array} \right. \qquad s=\overline{1,Q},  
\end{equation}
while the situations with~$rank(A_{pq})=Q-1$ feature both types of local behaviors. 
\end{remark}
 
\begin{lemma}
A sufficient condition for the occurence of a Dirac point~\cite{TWZ17}, i.e.~conical contact between dispersion surfaces, at~$(\bk^\ba,\omega_n^\ba)$ is given by  
\[ 
rank(A_{pq})=Q \quad \forall \hat\bk/\|\hat{\bk}\| \in \{\boldsymbol{\kappa}\in\mathbb{R}^d\!:\, \|\boldsymbol{\kappa}\|=1\}, \quad Q>1. 
\]
As a result, ``simple'' ($Q=2$) Dirac points of the scalar wave equation are not possible. 
\begin{proof} When~$Q=2$, the eigenvalues of~$A_{pq}$ are~$\pm\bth_{12}^{\mbox{\tiny{(0)}}}\!\cdot\hat{\bk}$, which vanish for~$\hat\bk/\|\hat{\bk}\|\perp\bth_{12}^{\mbox{\tiny{(0)}}}$. \end{proof}
\end{lemma}

\begin{remark}
The effective description given by Theorem~\ref{61} inherently depends upon the choice of the eigenfunction basis~$\tilde{\varphi}^{\ba}_{nq}$ ($q=\overline{1,Q}$), since any linear combination thereof is also an eigenfunction. To shed light on the problem, let~$T_{pq}$ be an invertible real-valued matrix such that 
\begin{equation} \label{transfo}
\widetilde{\tilde{\varphi}^{\ba}_{np}} ~=~ \sum_q T_{pq} \hh \tilde{\varphi}^{\ba}_{nq}, \qquad p=\overline{1,Q}
\end{equation}
represents another (possibly rescaled) basis that is orthogonal in~$L^2_\rho(Y_\ba)$. In this new frame, we have
\begin{equation} \label{modulitrans}
\widetilde{\rho_{p}^{\mbox{\tiny{(0)}}}} ~=~  \sum_r T_{pr} \hh\rho_{r}^{\mbox{\tiny{(0)}}}\hh T_{rp}^{\hh\mbox{\tiny{T}}}, \qquad
\widetilde{\bth_{pq}^{\mbox{\tiny{(0)}}}} ~=~ \sum_{r,s} T_{pr} \hh\bth_{rs}^{\mbox{\tiny{(0)}}}\hh T_{sq}^{\hh\mbox{\tiny{T}}}, \qquad 
\widetilde{\bmu_{pq}^{\mbox{\tiny{(0)}}}} ~=~ \sum_{r,s} T_{pr} \hh\bmu_{rs}^{\mbox{\tiny{(0)}}}\hh T_{sq}^{\hh\mbox{\tiny{T}}},
\end{equation}
due to~\eqref{ADmat} and~\eqref{Bmat}. As a result, effective models~\eqref{even5} and~\eqref{odd7A0} hold under the mappings  
\begin{multline} \label{transforule}
\ww_{jq} \:\to\: \widetilde{\ww_{jq}} = \sum_s T_{qs}^{\hh\mbox{\tiny{-T}}} \hh \ww_{js}, \qquad M_{pq}\:\to\: \widetilde{M_{pq}} = \sum_{r,s} T_{pr} M_{rs} T_{sq}^{\hh\mbox{\tiny{T}}}, \\ 
\langle e^{i\bk^\ba\cdot\bx}\rangle_{\ba}^{p\varphi} \:\to\: \langle e^{i\bk^\ba\cdot\bx}\rangle_{\ba}^{p\widetilde{\varphi}} = 
\sum_s T_{qs} \hh \langle e^{i\bk^\ba\cdot\bx}\rangle_{\ba}^{s\varphi},   \qquad  j\in\{0,1\}, \quad M\in\{A,B,D\}. 
\end{multline}
which clearly preserve the germane dispersion relationships. 
\end{remark}

The foregoing remark raises the question of a unique ``reference'' eigenfunction basis that facilitates comparison of results across studies. To help the discussion, we let
\begin{equation} \label{rho0}
D_{pq}^{\mbox{\tiny{-1/2}}} \;=\; \delta_{pq}\hh (\rho_{q}^{\mbox{\tiny{(0)}}})^{\mbox{\tiny{-1/2}}}, \qquad  \rho^{\mbox{\tiny{(0)}}} \;=\; \frac{1}{Q} \sum_q \rho^{\mbox{\tiny{(0)}}}_q,
\end{equation}
and we denote by $\lambda^{\mbox{\tiny{$A$}}}_p \in \mathbb{R}$ and $\lambda^{\mbox{\tiny{$B$}}}_p \in \mathbb{R}$ ($p=\overline{1,Q}$) the real eigenvalues of~$\hh  \rho^{\mbox{\tiny{(0)}}}\sum_{r,s}D_{pr}^{\mbox{\tiny{-1/2}}} A_{rs} D_{sq}^{\mbox{\tiny{-1/2}}}\hh$ and~$\hh  \rho^{\mbox{\tiny{(0)}}}\sum_{r,s}D_{pr}^{\mbox{\tiny{-1/2}}} B_{rs} D_{sq}^{\mbox{\tiny{-1/2}}}$, respectively. Note that~$\lambda^{\mbox{\tiny{$A$}}}_p$ are ordered in~$\pm$ pairs as in Remark~\ref{rem12}. 

\begin{lemma} \label{diagonal}
Let~$\hat\bk/\|\hat{\bk}\|$ be fixed, and consider the eigenfunction basis~$\widetilde{\tilde{\varphi}^{\ba}_{np}}$ given by~\eqref{transfo} with  
\begin{equation}\label{special}
T_{pq} ~=~ \sqrt{\rho^{\mbox{\tiny{(0)}}}} \sum_{s} R_{ps}^{\mbox{\tiny{T}}} \hh D_{sq}^{\mbox{\tiny{-1/2}}},  \qquad p,q=\overline{1,Q},  
\end{equation}
where~$R_{pq}$ is a real-valued orthogonal matrix that: (a)~converts $\rho^{\mbox{\tiny{(0)}}}\sum_{r,s}D_{pr}^{\mbox{\tiny{-1/2}}} A_{rs} D_{sq}^{\mbox{\tiny{-1/2}}}$ into a $2\!\times\!2$ block-diagonal form~$(\cdot)_{ot}$ via transformation~$\sum_{p,q} R_{op}^{\mbox{\tiny{T}}}\hh(\boldsymbol{\cdot})_{pq}\hh R_{qt}$ when solving~\eqref{even5}, and (b)~collects the eigenvectors of $\rho^{\mbox{\tiny{(0)}}}\sum_{r,s}D_{pr}^{\mbox{\tiny{-1/2}}} B_{rs} D_{sq}^{\mbox{\tiny{-1/2}}}$ when solving~\eqref{odd7A0}. Such basis is (i)~$\rho$-orthonormal in that~$(\rho\widetilde{\tilde{\varphi}^{\ba}_{np}}, \widetilde{\tilde{\varphi}^{\ba}_{nq}})_{\overline{Y}_{\!\ba}} = \delta_{pq}\hh \rho^{\mbox{\tiny{(0)}}}$ and (ii)~unique, i.e. independent of the choice of~$\tilde{\varphi}^{\ba}_{np}$, when~$M_{pq}$ itself ($M=A,B$) has no repeated eigenvalues. In this frame of reference, system~\eqref{even5} (resp. system~\eqref{odd7A0}) becomes $2\!\times\!2$ block-diagonal (resp. diagonal) in direction~$\hat\bk/\|\hat{\bk}\|$ as   
\begin{eqnarray} \label{uncoupledpde}
\widetilde{A_{pq}} &\!\!\!=\!\!\!& \widetilde{\bth_{pq}^{\mbox{\tiny{(0)}}}}\!\cdot\!(i\hat{\bk}) ~=~ \delta_{p(q+(-1)^p)}\hh (-1)^p i\hh \lambda^{\mbox{\tiny{$A$}}}_p, \notag \\ 
\widetilde{B_{pq}} &\!\!\!=\!\!\!& \widetilde{\bmu_{pq}^{\mbox{\tiny{(0)}}}}:(i\hat{\bk})^2 ~=~ \delta_{pq}\hh  \lambda^{\mbox{\tiny{$B$}}}_p, \qquad 
\widetilde{D_{pq}} ~=~ \delta_{pq}\hh \rho^{\mbox{\tiny{(0)}}} \qquad \text{(no summation)} 
\end{eqnarray} 
where~$\widetilde{A_{pq}}, \widetilde{B_{pq}}$ and~$\widetilde{D_{pq}}$ are invariant with respect to the choice of basis~$\tilde{\varphi}^{\ba}_{np}$. See Appendix~A for proof. 
\end{lemma}

With reference to~\eqref{uncoupledpde}, we note that the~$2\!\times\!2$ block-diagonal nature of~$\widetilde{A_{pq}}$ is perhaps not surprising, for in this case the germane dispersion branches appear in ``$\pm$'' pairs in terms of the incipient group velocity, see Remark~\ref{rem12}. We also note, omitting the details for brevity, that the analysis of Lemma~\ref{diagonal} also applies to ``reduced'' system~\eqref{odd9} arising when~$rank(A_{pq})=Q-1$.

\begin{remark}
At this point, it is unclear under which conditions there exists a common basis $\widetilde{\tilde{\varphi}^{\ba}_{np}}$ $(p=\overline{1,Q})$ that simultaneously \mbox{(block-)} diagonalizes~$M_{pq}$ ($M=A,B$) for all~$\hat\bk/\|\hat{\bk}\|$. For instance when Dirac cones are present $(M\!=\!A)$, there is no such basis; otherwise for each~$\hat{\bk}\perp\widetilde{\bth_{pq}^{\mbox{\tiny{(0)}}}}$ $(p,q=\overline{1,Q})$ the rank of~$\widetilde{A_{pq}}$ drops by at least two, and thus negates the Dirac cones. For other types of apexes (with~$M\!=\!B$) examined in Section~\ref{numerics} and elsewhere in the literature~\cite{Cra11}, on the other hand, such common basis does seem to exist. 
\end{remark}

\section{Nearby eigenvalues}\label{Near}

\noindent So far we considered the effective wave motion in a neighborhood of~$(\bk^\ba, \tilde{\lambda}^\ba_n)$, where the eigenvalue $\tilde{\lambda}^\ba_n$ is either simple, or repeated with multiplicity~$Q$. An implicit assumption in the underpinning analyses, however, is that the distance between $\tilde{\lambda}^\ba_n$ and its neighbors is~$O(1)$. In what follows, we relax this restriction by considering a \emph{cluster of $Q$ nearby eigenvalues} $\tilde{\lambda}^\ba_j$, $j=\overline{n\!-\!n_a,n\!+\!n_b}$ that are separated by~$O(\eps)$. For generality, we allow for the possibility that some of the eigenvalues $\tilde{\lambda}^\ba_j$ in this set may be equal, i.e. repeated. 

\begin{remark} To aid the analysis, we denote the featured cluster by~$\tilde{\lambda}^\ba_{j(q)}$, ($q=\overline{1,Q}$) where the mapping $j(q):[1,Q]\to [n-n_a,n+n_b]$ is one-to-one, albeit not necessarily monotonic. In the sequel, we also let $\phaq$ ($q=\overline{1,Q}$) be the eigenfunctions corresponding respectively to $\tilde{\lambda}^\ba_{j(q)}$. 
\end{remark}

\subsection{Eigenfunction basis} \label{ebasis2x}

In line with the above description, we consider a set of nearby eigenvalues $\tilde{\lambda}^\ba_{j(q)}$ given by 
\begin{equation}\label{cluster1}
\tilde{\lambda}^{\ba}_{j(q)} ~=~ \tilde{\lambda}^{\ba}_{n} + \eps\hh \gamma_q, \quad \gamma_q \lesseqgtr 0, \quad q=\overline{1,Q}.      
\end{equation}

\begin{assum}\label{repeat1} We allow for~$\tilde{\lambda}^{\ba}_{n}$ to have multiplicity \mbox{$1\leqslant Q'\!<Q$}, and we conveniently select the mapping~$j(q)$ so that $\gamma_q=0\,$ for $q=\overline{1,Q'}$. All other eigenvalues $\tilde{\lambda}^\ba_{j(q)}$ ($q>Q'$) are assumed to be simple. 
\end{assum}

From~\eqref{cluster1}, one finds that the eigenfunctions $\phaq \in H^1_p(Y)$ corresponding to~$\tilde{\lambda}^{\ba}_{j(q)}$  solve 
\begin{eqnarray} \label{DoubleEigeignpdex}
&&\hspace*{-5mm} -(\nabla\!+\!i\bk^\ba) \!\cdot\! \big(G (\nabla\!+\!i\bk^\ba) \phaq \big) ~=~ \tilde{\lambda}^{\ba}_n \hh\rho\hh \phaq + \eps\hh\gamma_q\hh\rho\hh \phaq  \quad \text{in}~~ Y, \\
&& \bnu\cdot G (\nabla\!+\!i\bk^\ba)\phaq |_{x_j=0} ~=\; -\bnu\cdot G (\nabla\!+\!i\bk^\ba)\phaq |_{x_j=\ell_j}. \label{Doubleeignbcx}
\end{eqnarray}
Following the previous analysis, we introduce the real-valued eigenfunctions~$\vphaq \in H^1_p(Y_\ba)$ that are: (i) orthogonal in~$L^2_\rho(Y_\ba)$, (ii) related to~$\vphaq$ via~\eqref{Eigensystem-apex2}, and (iii) solve 
\begin{eqnarray} \label{Eigensystem-apex412x}
&& \hspace*{-5mm} -\nabla\!\cdot\! \big(G\nabla\vphaq\big) ~=~ \tilde{\lambda}^{\ba}_{n} \hh \rho \, \vphaq + \eps\hh\gamma_q\hh\rho\hh \vphaq \quad \text{in}~~ Y_\ba, \\*[-0.8mm]
&& \begin{array}{rcl} \label{Eigensystem-apex432x}
a_j \vphaq|_{x_j=0} &\!\!\!=\!\!\!& -a_j \vphaq|_{x_j=\ell_j},  \\*[1.2mm]
\bnu\cdot G\nabla\vphaq|_{x_j=0} &\!\!\!=\!\!\!& -\bnu\cdot G\nabla\vphaq|_{x_j=(1+a_j)\ell_j}.
\end{array}
\end{eqnarray}
Inherently, \eqref{DoubleEigeignpdex}--\eqref{Eigensystem-apex432x} suggest use of the ``averaging'' operators~\eqref{effectiveab2} and function space~\eqref{hilbertap}.

\subsection{Leading-order approximation} \label{Loax}

We again refer to~\eqref{ComparisonAnsatzEqna}--\eqref{ComparisonAnsatzEqnBCNa}, and we start from the asymptotic expansion~\eqref{wexpa} of~$\tw$. From the results in Section~\ref{Bri}\ref{loaa}, one finds that the $O(\eps^{-2})$ contribution to~\eqref{ComparisonAnsatzEqna} is given by~\eqref{w0ea}. Thanks to~\eqref{Eigensystem-apex412x}--\eqref{Eigensystem-apex432x}, one finds that~\eqref{w0ea} is satisfied by~\eqref{DoubleEigw0sol} up to the $O(\eps)$ residual 
\begin{equation}\label{resid0}
\eps\hh  \rho(\bx) \sum_q \gamma_q \hh \ww_{0q} \hh \vphaq(\bx).  
\end{equation}
This remainder inherently carries over to the $O(\eps^{-1})$ statement of~\eqref{ComparisonAnsatzEqna}, which now reads  
\begin{eqnarray}\label{near2}
-\tilde{\lambda}^{\ba}_n\hh\rho\hh\tw_1 - \sigma\breve{\omega}^2\!\rho\hh\tw_0
- \nabla \!\cdot\! \big(G(\nabla\tw_1 \!+ i\hat{\bk}\hh\tw_0)\big) - i \hat{\bk}\hh \!\cdot\! \big( G \nabla \tw_0\big) - \rho \sum_q \gamma_q \hh \ww_{0q} \hh\vphaq ~=~ 
0 \quad \text{in}~~ Y_{\ba}.~
\end{eqnarray}
On integrating~$\langle\eqref{near2}\rangle_\ba^{p\varphi}$ by parts, one can show as in the proof of Lemma~\ref{lem3} that 
\begin{equation} \label{GEPgam}
\sum_{q} A_{pq}^\gamma \hh\ww_{0q} \,-\, \tau \sum_{q} D_{pq} \hh\ww_{0q} ~=~0, \qquad p=\overline{1,Q} 
\end{equation}
up to respective $O(\eps)$ residuals
\begin{equation}\label{resid1eff}
\eps \gamma_p \hh\langle\rho\hh\tw_1\rangle_\ba^{p\varphi},  \qquad p=\overline{1,Q},   
\end{equation}
where 
\begin{equation} \label{ADmatgam}
A_{pq}^\gamma ~=~ \bth_{pq}^{\mbox{\tiny{(0)}}} \!\cdot\! i\hat{\bk} + \gamma_{q} \hh D_{pq}, \qquad D_{pq} ~=~ \delta_{pq} \hh \rho_{q}^{\mbox{\tiny{(0)}}}, \qquad \tau  = -\sigma \breve{\omega}^2,  
\end{equation} 
with~$\rho_{p}^{\mbox{\tiny{(0)}}}$ and~$\bth_{pq}^{\mbox{\tiny{(0)}}}$ given by~\eqref{pqmoduli}. As before, we seek (when possible) a non-trivial solution to~\eqref{GEPgam}. 

\begin{remark}\label{lem7gam}
Since $A_{pq}^\gamma$ is Hermitian and~$D_{pq}$ is diagonal positive definite, the eigenvalues~$\tau$ of~\eqref{GEPgam} are real, and the eigenvectors corresponding to distinct eigenvalues are orthogonal.  
\end{remark}

\subsubsection{Effective model for full-rank $A_{pq}^\gamma$} \label{Apqnot0x}

When $rank(A_{pq}^\gamma)=Q$, generalized eigenvalue problem~\eqref{GEPgam} has \emph{no zero eigenvalues}. Accordingly, the scaling framework~\eqref{FW-FF} with~\mbox{$\breve{\omega}=0$} does not cater for a non-trivial solution of~\eqref{GEPgam}; we thus keep~$\breve{\omega}\neq0$ and set~$\ww_{0q}\!=0$ $(q=\overline{1,Q})$ by which the $O(\eps^{-1})$ statement~\eqref{near2} reduces to~\eqref{even1}. Thanks to~\eqref{Eigensystem-apex412x}--\eqref{Eigensystem-apex432x}, one can take~$\tw_1(\bx)$ as in~\eqref{even2}, which satisfies~\eqref{even1} up to the $O(\eps)$ residual 
\begin{equation}\label{resid1}
\eps\hh \rho(\bx) \sum_q \gamma_q \hh \ww_{1q} \hh \vphaq(\bx),   \qquad 
\end{equation}
that carries over to the~$O(1)$ statement of~\eqref{ComparisonAnsatzEqna}, namely 
\begin{multline}\label{even3gam}
-\tilde{\lambda}^{\ba}_n\hh\rho\hh\tw_2 - \nabla \!\cdot\! \big(G (\nabla\tw_2 
\!+ i\hat{\bk}\hh\tw_1)\big) - i \hat{\bk}\hh \!\cdot\! \big( G \nabla \tw_1\big) - \sigma \breve{\omega}^2\rho\hh\tw_1 \\
- \rho\hh \sum_q \gamma_q \hh \ww_{1q} \hh\vphaq \:=\: e^{i\bk^\ba\cdot\bx} \quad \text{in}~~ Y_{\ba}.
\end{multline}
To compute constants~$\ww_{1q}\in\mathbb{C}$ in~\eqref{even2}, we integrate~$\langle\eqref{even3gam}\rangle_{\ba}^{p\varphi}$ $\hh (p\in\overline{1,Q})\hh$ by parts which yields 
\begin{eqnarray} \label{even5gam}
-\sum_{q} A_{pq}^\gamma \hh \ww_{1q} \,-\, \sigma \breve{\omega}^2 \sum_{q} D_{pq}\hh \ww_{1q} ~=~ 
\langle e^{i\bk^\ba\cdot\bx}\rangle_{\ba}^{p\varphi}, \qquad p=\overline{1,Q} 
\end{eqnarray}
up to $O(\eps)$ residuals, where~$A_{pq}^\gamma$ and~$D_{pq}$ are given by~\eqref{ADmatgam}, and $\sigma \breve\omega^2 = \eps^{-1}(\omega^2\!-\tilde{\lambda}^{\ba}_{n})$ according to~\eqref{FW-FF}. By Remark~\ref{lem7gam}, we find that the generalized eigenvalue problem behind~\eqref{even5gam} has~$Q$ real eigenvalues, hereon denoted by~$\tau_r$ ($r=\overline{1,Q}$), and we assume~$-\sigma\breve{\omega}^2 \neq \tau_r$ when solving~\eqref{even5gam}. Since~$A_{pq}^\gamma$ is not antisymmetric, however, the eigenvalues~$\tau_r$ of~\eqref{even5gam} will generally not appear in ``$\pm$'' pairs as in the case of repeated eigenvalues, see Section~\ref{Rep}\ref{Loa}\ref{Apqnot0}.  

\begin{remark} The case of trivial~$A_{pq}^\gamma=A_{pq}+\gamma_q D_{ps}$ is not germane to the discussion for it would require that each $A_{pq}\in\mathbb{I}$ and $\gamma_q D_{pq}\in\mathbb{R}$ vanish identically -- a premise that yields $\gamma_p=0$ ($p=\overline{1,Q}$) and thus precludes the existence of nearby eigenvalues.  
\end{remark}

\subsubsection{Effective model for $rank(A_{pq}^\gamma)=Q\!-\!1$} \label{ApqQm1}

\noindent For clarity we focus on situations where, upon suitable ``rotation'' \eqref{transfo} of the eigenfunction basis~$\vphaq$ -- as permitted by the multiplicity (if any) of eigenvalue $\tilde{\lambda}^\ba_n=\tilde{\lambda}^\ba_{j(q)}$ ($q=\overline{1,Q'}$), one has     
\begin{equation} \label{A000gam}
A_{1p}^\gamma=A_{p1}^\gamma=0, \quad p = \overline{1,Q}, 
\end{equation}
cf.~\eqref{A000}. Recalling~\eqref{ADmatgam}, we note that~\eqref{A000gam} is consistent with the earlier premise that~$\gamma_1=0$, see Assumption~\ref{repeat1}. 

\begin{lemma}\label{Apqgamma}
Let $rank(A_{pq}^\gamma)=Q\!-\!1$, and let the direction~$\hat{\bk}/\|\hat{\bk}\|$ be fixed. If~\eqref{A000gam} holds for some \mbox{$\kappa=\|\hat{\bk}\|>0$}, then it holds for all $\kappa>0$. 
\begin{proof}
From~\eqref{ADmatgam}, we see that 
\[
A_{pq}^\gamma ~=~ \gamma_q\hh D_{pq} \,+\, (\bth_{pq}^{\mbox{\tiny{(0)}}}\!\cdot \hat{\bk}/\kappa) \hh i \kappa, \quad p,q\in\overline{1,Q}
\]
where $\gamma_q\hh D_{pq}\in\mathbb{R}$ and $(\bth_{pq}^{\mbox{\tiny{(0)}}}\!\cdot \hat{\bk}/\kappa)\in\mathbb{R}$ are independent of~$\kappa$. If~\eqref{A000gam} holds for some $\kappa>0$, then each of these two terms must vanish identically for  $(p,q)\to(1,p)$ and $(p,q)\to(p,1)$.  
\end{proof}
\end{lemma}

By virtue of~\eqref{A000gam}, \eqref{GEPgam} becomes identity if: (a) $\hh\tau = -\sigma\breve{\omega}^2=0$, and (b) $\hh w_{0q} = \text{w}_0\hh \delta_{q1}$ i.e. 
\begin{equation} \label{new2x}
\tw_{0}(\bx) ~=~ \text{w}_0 \hh \vphao(\bx), 
\end{equation}
for some constant~$\text{w}_0\in\mathbb{C}$. On recalling that~$\gamma_1=0$, this result yields $\gamma_q\hh w_{0q}=0$ ($q=\overline{1,Q}$); accordingly the $O(\eps^{-1})$ field equation~\eqref{near2} is solved by 
\begin{equation}\label{w1solacx}
\tw_1(\bx) ~=~ \text{w}_0 \bchi_{1}^{\mbox{\tiny{(1)}}}(\bx)\!\cdot i\hat{\bk}  \;+\; \sum_q \ww_{1q}\hh\vphaq(\bx), 
\end{equation}
up to the~$O(\eps)$ residual~\eqref{resid1}. Here~$\ww_{1q}\!\in\mathbb{C}$ are constants, and $\bchi_1^{\mbox{\tiny{(1)}}}\in (H_{p0}^{1\ba}(Y_\ba))^d$ uniquely solves  
\begin{eqnarray} \label{Comparisonchi1a1xx}
\hspace*{-7mm} && \tilde{\lambda}^{\ba}_{n}\hh\rho\hh\bchi_1^{\mbox{\tiny{(1)}}} + 
\nabla \!\cdot\! \big(G(\nabla\bchi_1^{\mbox{\tiny{(1)}}}\! + \bI \vphao)\big) + G \nabla \vphao 
- \sum_r\frac{1}{\rho_{r}^{\mbox{\tiny{(0)}}}}\hh\bth_{r1}^{\mbox{\tiny{(0)}}}\hh\rho\hh \vphar ~=~0 \quad \text{in}~~ Y_{\ba}, \qquad \\*[-0.5mm] \notag
\hspace*{-7mm} && \quad\qquad \bnu \cdot G(\nabla\bchi_1^{\mbox{\tiny{(1)}}}\! + \bI \vphao)|_{x_j=0} ~=\: 
-\bnu \cdot G(\nabla\bchi_1^{\mbox{\tiny{(1)}}}\! + \bI \vphao)|_{x_j=(1+a_j)\ell_j}.   
\end{eqnarray}
With regard to~\eqref{Comparisonchi1a1xx}, we in particular note that $\tilde{\lambda}^{\ba}_{j(1)}=\tilde{\lambda}^{\ba}_{n}$ because~$\gamma_1=0$. In this setting the~$O(1)$ statement of~\eqref{ComparisonAnsatzEqna}, which accounts for residual~\eqref{resid1}, reads 
\begin{multline}\label{even3gamb}
-\tilde{\lambda}^{\ba}_n\hh\rho\hh\tw_2 - \nabla \!\cdot\! \big(G (\nabla\tw_2 
\!+ i\hat{\bk}\hh\tw_1)\big) - i \hat{\bk}\hh \!\cdot\! \big(G(\nabla\tw_1 + i \hat{\bk}\hh \tw_0)\big) - \sigma \hat{\omega}^2\rho\hh\tw_0 \\ 
- \rho \sum_q \gamma_q \hh \ww_{1q} \hh\vphaq \:=\: e^{i\bk^\ba\cdot\bx} \quad \text{in}~~ Y_\ba.  
\end{multline}
On integrating $\langle\eqref{even3gamb}\rangle_\ba^{1\varphi}$ by parts, we obtain 
\begin{equation} \label{odd7gam}
- \big(\bmu_{1\nes 1}^{\mbox{\tiny{(0)}}}\!:\!(i\hat{\bk})^2 + 
\sigma\hat{\omega}^2\rho_{1}^{\mbox{\tiny{(0)}}}\big) \text{w}_0 ~=\, 
 \langle e^{i\bk^\ba\cdot\bx}\rangle_{\ba}^{1\varphi}, 
\end{equation}
where~$\rho_{1}^{\mbox{\tiny{(0)}}}$ and~$\bmu_{1\nes 1}^{\mbox{\tiny{(0)}}}$ are given respectively by~\eqref{pqmoduli} and~\eqref{Bmat}. Note that~\eqref{odd7gam} is formally the same as~\eqref{odd7}, catering for the case of repeated eigenvalues. 

To identify the remaining~$Q\!-\!1$ leading-order projections in terms of~$w_{1q}$ $(q=\overline{2,Q})$, we apply the analysis from Section~\ref{Rep}\ref{Loa}\ref{Apqneg1} to~\eqref{even3gamb} with~$\breve{\omega}\neq 0$ to obtain the linear system 
\begin{multline} \label{odd9gam}
-\sum_{q>1} A_{pq}^\gamma\hh \ww_{1q} \,-\, \sigma \breve{\omega}^2 \sum_{q>1} D_{pq}\hh \ww_{1q} ~=~ 
\langle e^{i\bk^\ba\cdot\bx}\rangle_{\ba}^{p\varphi} \:+\: 
\text{w}_0 \hh \big(\bmu_{p1}^{\mbox{\tiny{(0)}}}\!:\!(i\hat{\bk})^2 + \sigma\breve{\omega}^2 \brho_{p1}^{\mbox{\tiny{(1)}}} \!\cdot\!i\hat{\bk}\big), \\*[-2.7mm] 
\brho_{p1}^{\mbox{\tiny{(1)}}} = \langle\rho\hh\bchi_{1}^{\mbox{\tiny{(1)}}}\rangle_{\ba}^{p\varphi}, \qquad p=\overline{2,Q}.
\end{multline} \vspace*{-6mm}

\begin{remark}
We observe a remarkable similarity between~\eqref{odd9gam} and its repeated-eigenvalues counterpart~\eqref{odd9}, where the only formal difference is the appearance of a diagonally-perturbed matrix, $A_{pq}^\gamma=A_{pq}+\gamma_q D_{pq}$, in lieu of~$A_{pq}$. The same comment applies to~\eqref{even5gam} and~\eqref{even5}. As will be seen shortly, the perturbation contained in $A_{pq}^\gamma$ has the effect of ``blunting'' the conical dispersion surfaces. 
\end{remark} \vspace*{-2mm}

Analogous to the case of repeated eigenvalues, the ``full'' effective model when~$rank(A_{pq}^\gamma)=Q\!-\!1$ is given by~$\text{w}_{0}$ solving~\eqref{odd7gam}, and~$w_{1q}$ ($q=\overline{2,Q}$) solving~\eqref{odd9gam}. In terms of the implied dispersion within the cluster, we further see that~\eqref{odd7gam}--\eqref{odd9gam} describe a single parabola accompanied by \mbox{$Q\!-\!1$} ``cones''. Note that when computing the ``conical'' dispersion due to~\eqref{odd9gam}, one must set~$\text{w}_0=0\hh$ because (in the absence of a source term) its nontrivial values  require by~\eqref{odd7gam} that $\omega^2\!-\tilde{\lambda}^{\ba}_n=O(\eps^2)$, whereas~\eqref{odd9gam} assumes $\hh\omega^2\!-\tilde{\lambda}^{\ba}_n=O(\eps)$. 


\section{Numerical results} \label{numerics}

\noindent We next illustrate the analysis by studying the effective wave motion in one-dimensional polyatomic chains and chessboard-like periodic continua in~$\mathbb{R}^2$. 

\begin{remark}
Given the lack of a priori knowledge as to when precisely the nearby eigenvalues model (Section~\ref{Near}) becomes relevant, we first restrict our asymptotic description to the (second-order) simple eigenvalue model and the repeated eigenvalue model described respectively Section~\ref{Bri} and~Section~\ref{Rep}. Later on, we exercise the results from Section~\ref{Near} to extend the range of validity of local approximation as called upon by the comparison with numerical simulations. 
\end{remark} 

\subsection{FW-FF homogenization of polyatomic chains} \label{example1}
 
\begin{figure}[h!]\vspace*{-2mm}
\centering
\centering{\includegraphics[width=0.9\linewidth]{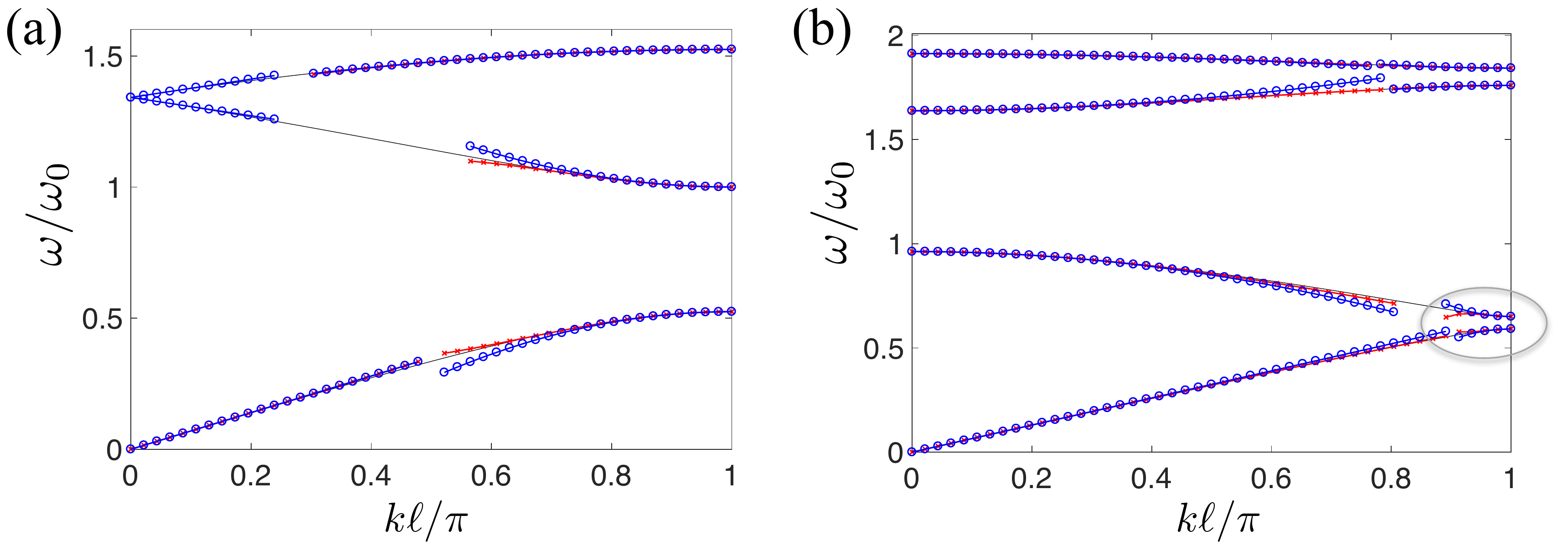}}\vspace*{-4mm}
\caption{Dispersion in polyatomic chains: (a) triatomic chain with $\boldsymbol{m}=(1,2.5,1)$ and $\boldsymbol{c}=(1,1,0.4)$, and (b) tetratomic chain with $\boldsymbol{m}=(1,2.0237,1.3539,1.6625)$ and $\boldsymbol{c}=(1,0.4257,1.0521,0.3731)$. Here $\omega_{0}=\sqrt{c_1/m_1}$; solid line plots the exact relationship, and ``o'' (resp.~``x'') marker tracks the leading- (resp.~second-) order approximation.} \vspace*{-7mm}
\label{poly_chain}
\end{figure}

\noindent Here, we consider a periodic chain of alternating masses and springs whose unit cell of length~$\ell$ contains~$N$ equidistant masses~$m_{j}$ ($j=\overline{1,N}$). We let~$c_j$ be the stiffness of spring connecting~$m_j$ and~$m_{j+1}$. Application of the foregoing theory to polyatomic chains is summarized in Appendix~B (electronic supplementary material); see also~\cite{CKP10b} for a comprehensive leading-order analysis of lattice structures. Using the short-hand notation $\boldsymbol{m}=(m_1,\ldots m_N)$ and $\boldsymbol{c}=(c_1,\ldots c_N)$, Fig.~\ref{poly_chain} plots the FF-FW approximation of the dispersion relationships for two example polyatomic chains with~$N=3,4$. As can be seen from the display, the second-order model outperforms its leading-order counterpart and covers (for $N=3$) with reasonable accuracy at least 50\% of the first Brillouin zone. From panel~(a), it is further seen that a one-dimensional, discrete variant of the ``full-rank~$A_{pq}$'' model (Section~\ref{Rep}\ref{Loa}\ref{Apqnot0}) effectively describes the occurence of a double eigenvalue in a triatomic chain at~$k\ell=0$. The circled apex region in panel~(b), however, exposes shortcomings of the simple eigenvalue model (Section~\ref{Bri}) in situations involving nearby dispersion branches. One may further  note that the simple eigenvalue model describes well the third and the fourth branch of the tetratomic chain near~$k\ell=\pi$, despite the fact that their distance is \emph{commensurate} with that between the first and the second branch inside the circled region. This suggests that the relevance of the nearby eigenvalue model  (Section~\ref{Near}) is controlled by a parameter that entails both (i) distances between the neighboring eigenvalues, and (ii) incipient curvatures of the respective dispersion curves.

\subsection{FW-FF homogenization of a chessboard-like medium} \label{example2}

\begin{figure}[h!]\vspace*{-3mm}
\centering
\centering{\includegraphics[width=0.82\linewidth]{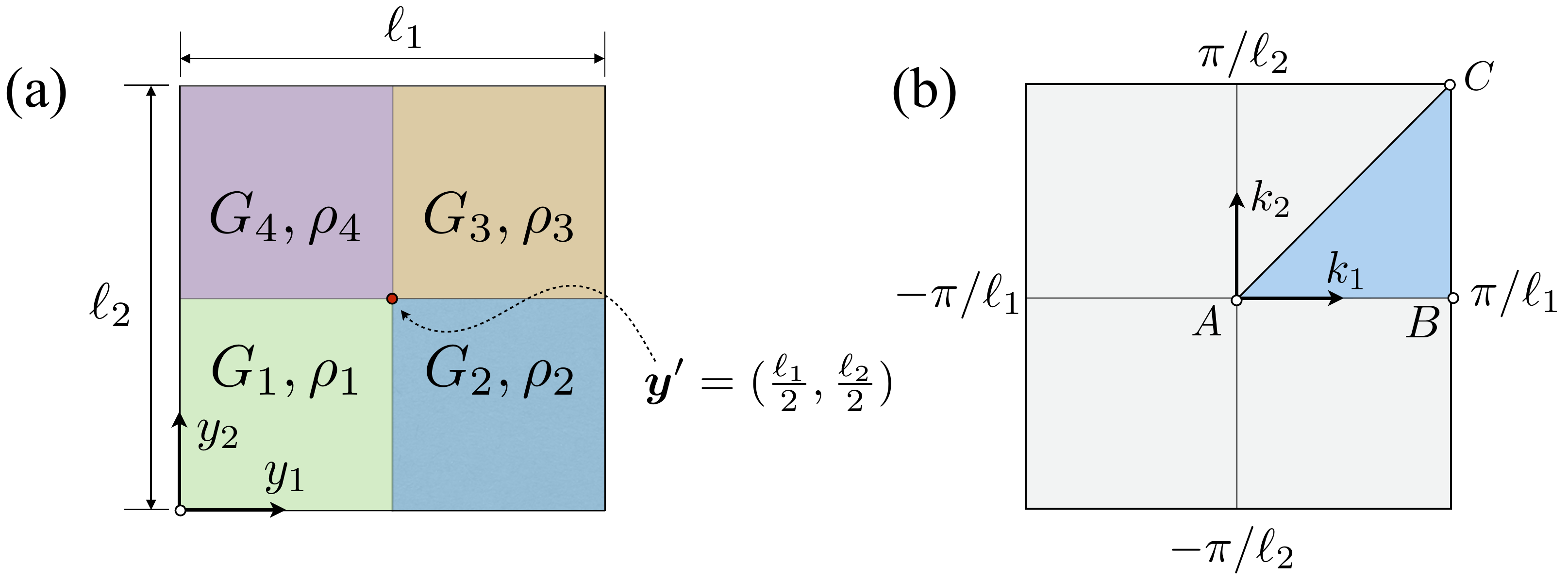}}\vspace*{-3mm}
\caption{Chessboard-like medium: (a) unit cell of periodicity~$Y$, and (b) first Brillouin zone~$\mathcal{B}$.} \vspace*{-2mm}
\label{chess_uc_bz}
\end{figure}

\noindent Consider a chessboard-like periodic medium whose unit cell~$Y\!\in\mathbb{R}^2\!$ and the first Brillouin zone~$\mathcal{B}$ are shown in Fig.~\ref{chess_uc_bz}. For brevity of notation, we let $\bG = (G_1,G_2,G_3,G_4)$ and $\br = (\rho_1,\rho_2,\rho_3,\rho_4)$.

\subsubsection{Green's function near the edge of a band gap} \label{example2green}

\begin{figure}[h!] \vspace*{-2mm}
\centering{\includegraphics[width=1.0\linewidth]{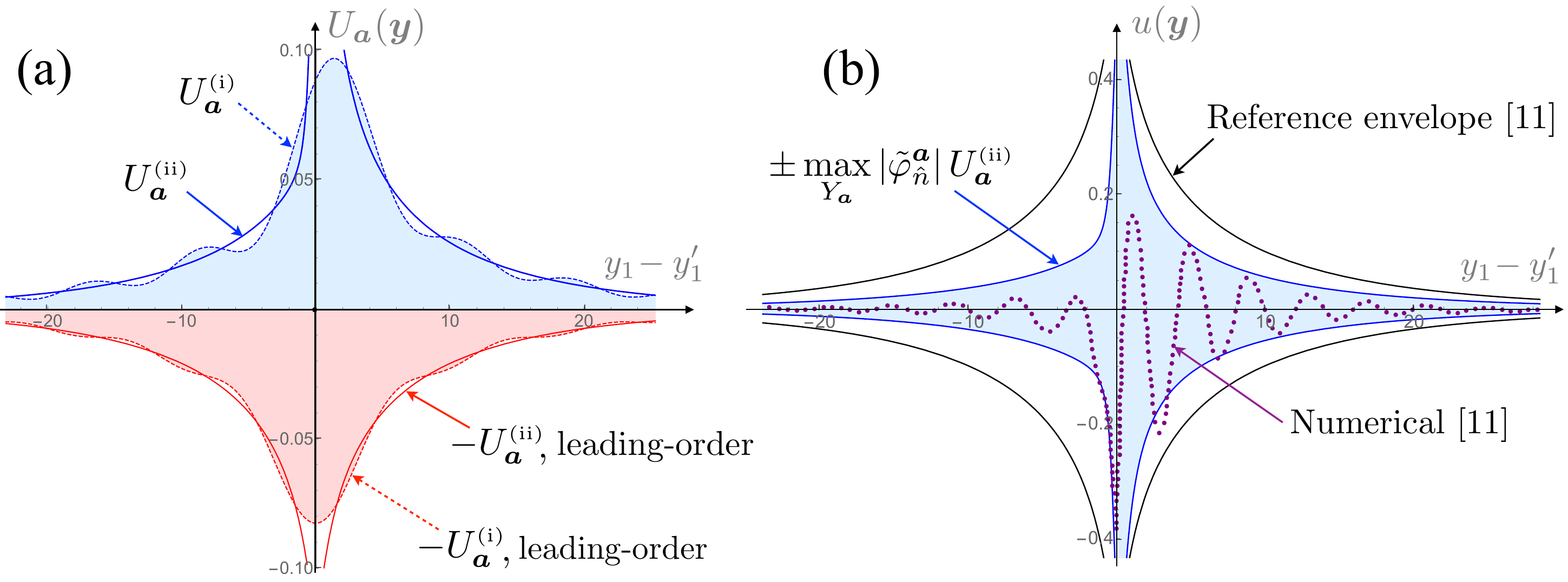}}\vspace*{-2mm}
\caption{Response of a chessboard-like medium with~$\ell_1\!=\!\ell_2\!=\!2$, $\hh\bG=(1,1,1,1)$ and~$\br=(1,101,201,101)$ along line~$y_2-y_2'=0.5$ due to the point source shown in Fig.~\ref{chess_uc_bz}(a) at frequency $\omega=0.1720$: (a) asymptotic fields~$U_\ba^{\mbox{\tiny{(i)}}}$ and~$U_\ba^{\mbox{\tiny{(ii)}}}$ solving~\eqref{gfun8} together with their leading-order (monopole) counterparts, and (b) numerical simulation of the problem~\cite{Cra11} versus ``trend'' envelope~$\pm\tfrac{1}{2\pi}K_1(\alpha\|\by\|)$ \cite{Cra11} and envelope~$\pm\max|\tilde{\varphi}_{\hat{n}}^\ba|\hh U_\ba^{\mbox{\tiny{(ii)}}}$ due to~\eqref{gfun3} and~\eqref{gfun8}.}\label{figog} \vspace*{-2mm}
\end{figure}

\noindent We first consider an elastic ``chessboard'' with~$\ell_1\!=\!\ell_2\!=\!2$, $\bG\!=\!(1,1,1,1)$,  and~$\br\!=\!(1,101,201,101)$ as in~\cite{Cra11}, and we seek an effective description of its response due to a point source acting at frequency $\omega=0.1720$ that is just inside the first band gap, see Fig.~\ref{figo1}(b). In this case, the results from Section~\ref{Bri}\ref{space} apply with~$\hat{n}=1$, $\ba=(1,1)$, and~$\epsilon^2=\omega^2-\tilde{\lambda}^{\ba}_1=1.374\times10^{-4}$. Placing the point source, $\delta(\by-\by')$, as in Fig.~\ref{chess_uc_bz}(a), we find that~$\vphah(\by')=0.9161$ and~$\bchi^{\mbox{\tiny{(1)}}}(\by')=(-1.704,-1.704)$. For such loading configuration, Fig.~\ref{figog}(a) compares the respective solutions~$U_\ba^{\mbox{\tiny{(i)}}}$ and~$U_\ba^{\mbox{\tiny{(ii)}}}$ of~\eqref{gfun6} with~$F(\by-\by')$ given by \eqref{gfun7} along line~$y_2-y_2'=\tfrac{1}{4}\ell_2$. Here $U_\ba^{\mbox{\tiny{(ii)}}}$ is given by~\eqref{gfun8}, while~$U_\ba^{\mbox{\tiny{(i)}}}$ is evaluated by integrating~\eqref{gfun4} numerically over~$\mathcal{B}_\ba$. For completeness, also included in the display are the \emph{leading-order} components of~$U_\ba^{\mbox{\tiny{(i)}}}$ and~$U_\ba^{\mbox{\tiny{(ii)}}}$, obtained by discarding the ``dipole'' contributions on the right hand-sides of~\eqref{gfun4} and~\eqref{gfun8}, i.e. those terms multiplied by~$\bchi^{\mbox{\tiny{(1)}}}(\by')$. We first note that approximation~(i) and approximation~(ii) give similar result away from~$y_1-y_1'=0$. The difference between~$U_\ba^{\mbox{\tiny{(i)}}}$ and~$U_\ba^{\mbox{\tiny{(ii)}}}$ on the one hand and their leading-order counterparts on the other is, however, striking in that the former bring about a \emph{directivity of wave motion} that is absent from the leading-order approximation. This is highlighted in Fig.~\ref{figog}(b), which compares the effective solution envelope~$\pm\max|\tilde{\varphi}_{\hat{n}}^\ba|\hh U_\ba^{\mbox{\tiny{(ii)}}}(\by)$ due to~\eqref{gfun3} and~\eqref{gfun8} with the corresponding numerical simulation~\cite{Cra11}. It is apparent that $\pm\max|\tilde{\varphi}_{\hat{n}}^\ba|\hh U_\ba^{\mbox{\tiny{(ii)}}}(\by)$ provides notably better effective description of the numerical simulation than the reference envelope~$\tfrac{1}{2\pi}K_1(\alpha\|\by\|)$ \cite{Cra11}, cf.~\eqref{gfun8}, in terms of both magnitude and directivity of the wave motion. For completeness the specifics of numerical implementation, including a tabulated comparison of the effective coefficients~$\rho^{\mbox{\tiny{(0)}}}$ and~$\bmu^{\mbox{\tiny{(0)}}}$ with those in ~\cite{Cra11}, is provided in Appendix~C (electronic supplementary material).

\subsubsection{Dispersion in a medium with variable shear modulus} \label{example2modulus}

\begin{figure}[h!] 
\centering{\includegraphics[width=0.88\linewidth]{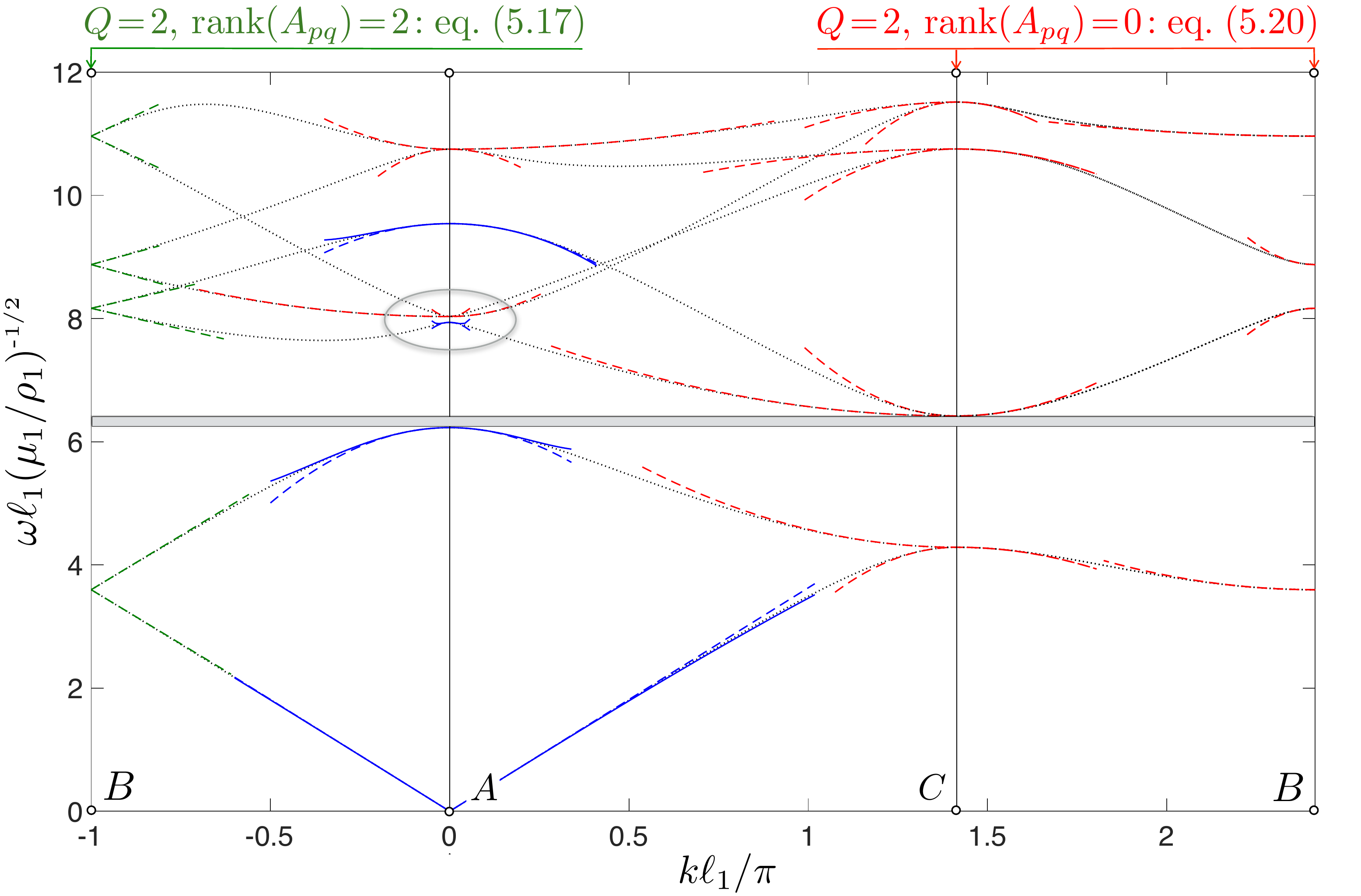}}\vspace*{-3.2mm}
\caption{Dispersion in a chessboard-like medium with~$\ell_1\!=\ell_2\!=1$, $\bG=(1,4,1,4)$, and~$\br=(1,2,1,2)$. Dotted lines plot the ``exact'' relationships computed numerically by NGSolve~\cite{NGSolve}, and dashed (resp.~continuous) lines track the leading- (resp.~second-) order approximations near apexes~$A$, $B$ and~$C$.} \vspace*{-3mm}
\label{chess_disper} \vspace*{-2mm}
\end{figure}

\noindent We next consider the dispersion in a chesboard-like medium with~$\ell_1=\ell_2=1$, $\bG=(1,4,1,4)$, and~$\br=(1,2,1,2)$. This configuration, examined in~\cite{MG17} from the LW-LF standpoint, differs from the previous example for it features a \emph{variable shear modulus}. Using the results from Section~\ref{Bri} and Section~\ref{Rep}, asymptotic approximations of the first eight dispersion branches computed and shown in Fig.~\ref{chess_disper} along the Brillouin contour~$ABC$, see Fig.~\ref{chess_uc_bz}(b). The leading-order approximations are plotted for all branches and all apex points, while their second-order counterparts are shown only for those pairs~$(\bk^\ba,\omega^\ba_n)$ featuring isolated eigenvalues~$\tilde{\lambda}^\ba_n$. As can be seen from the display, the second-order model brings about notable improvement in the asymptotic description of the first, second, and sixth branch near apex~$A$; however this observation does not apply to the third branch, see the circled apex region, due to the \emph{proximity of neighboring eigenvalues} (this issue will be addressed shortly). From Fig.~\ref{chess_disper}, one may also note that apex~$B$ features distinct asymptotic behaviors in directions~$BA$ and~$BC$. Indeed, in direction~$BA$ the local approximation is \emph{linear} ($\text{rank}(A_{pq})=2$) and the use is made of the generalized eigenvalue problem stemming from~\eqref{even5} to determine the germane slopes, see also Remark~\ref{dirac}. In direction~$BC$, on the other hand, the local variation is \emph{quadratic} ($\text{rank}(A_{pq})=0$) and the germane curvatures are computed from the generalized eigenvalue problem due to~\eqref{odd7A0}. Finally, we observe that the local behaviors of the first eight branches at apex~$C$ are uniformly described by the ``trivial~$A_{pq}$'' model~\eqref{odd7A0}. 

\subsection{Clusters of nearby dispersion branches} \label{example2clusters}

\begin{figure}[h!] \vspace*{-4.5mm}
\centering{\includegraphics[width=0.99\linewidth]{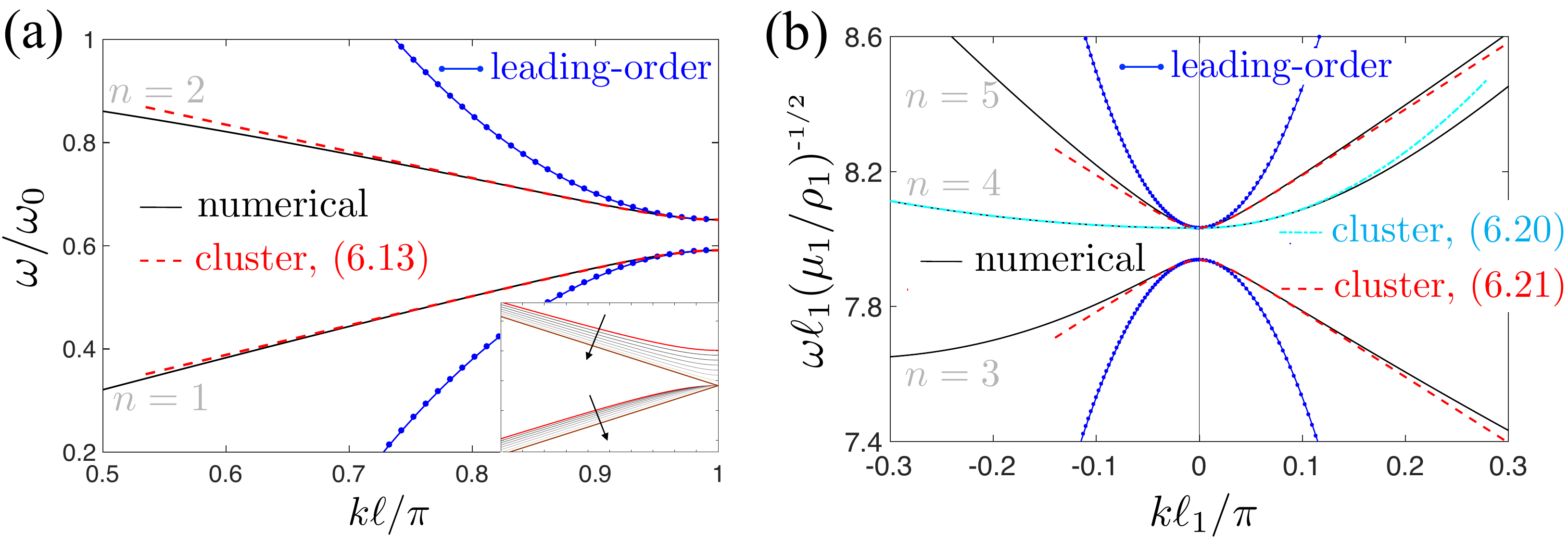}}\vspace*{-3mm}
\caption{Asymptotics of closely spaced dispersion branches: (a) two nearby eigenvalues at apex $k\ell=1$ in the tetratomic chain from Fig.~\ref{poly_chain}, and (b) three nearby eigenvalues at apex~$A$ in the chessboard medium from Fig.~\ref{chess_disper}.} 
\label{nearby} \vspace*{-4mm}
\end{figure}

To conclude the study, we apply the results from Section~\ref{Near} to examine the asymptotics of wave dispersion in situations featuring tightly spaced dispersion branches, see the circled regions in Fig.~\ref{poly_chain} and Fig.~\ref{chess_disper}. Considering first the tetratomic chain example with~$Q=2$ nearby eigenvalues at~$k\ell=\pi$, we find that $rank(A_{pq}^\gamma)=2$ whereby the ``\emph{twin cones}'' model~\eqref{even5gam} applies. On the other hand, for the cluster of~$Q=3$ nearby eigenvalues at apex~$A$ of the chessboard medium, one finds that $rank(A_{pq}^\gamma)=2$ and~\eqref{A000gam} holds, which implies the ``\emph{parabola with cones}'' model~\eqref{odd7gam}--\eqref{odd9gam}. Fig.~\ref{nearby} illustrates the performance of both models, from which one observes a major improvement in the asymptotic description of closely spaced dispersion branches. In particular, \eqref{even5gam} and~\eqref{odd9gam} are seen to be effective in representing  locally the \emph{``blunted cone''} shape of the respective dispersion curves. For further insight, the insert in Fig.~\ref{nearby}(a) depicts the evolution of dispersion surfaces due to~\eqref{even5gam} when the diagonal perturbation, $\gamma_{pq}D_{pq}$, of $A_{pq}^\gamma$ is artificially reduced to zero. Concerning Fig.~\ref{nearby}(b), one may also note from~\eqref{odd7gam}--\eqref{odd9gam} that the fourth branch (``parabola'') and the fifth branch (``cone'') \emph{do not interact} despite touching each other, while the third and the fifth branch (``cones'') interact despite separation. For completeness, details of calculation for both cases are provided in Appendix~D (electronic supplementary material). 

\section{Summary} \label{summary}

\noindent In this work, we establish a general framework for dynamic homogenization of wave motion in~$\mathbb{R}^d$, $d\geqslant 1$, at finite wavenumbers and frequencies. To cater for the nature of the problem, we pursue the plane-wave expansion approach and we adopt the projection of a Bloch wave onto an eigenfunction for the unit cell of periodicity -- at a fixed wavenumber and frequency -- as effective descriptor of wave motion. In this way we obtain a homogenized field equation (including the source term), for an arbitrary dispersion branch, near corners of the ``wavenumber quadrants'' comprising the first Brillouin zone. We consider the situations of both (i)~isolated, (ii)~repeated, and (iii)~nearby eigenvalues. The second-order analysis of isolated eigenvalues demonstrates that the leading correction of the effective differential operator is~$O(\eps^2)$. In contrast, the effective source term is found to undergo an~$O(\epsilon)$ correction, which may be important for the development of homogenized Green's functions for periodic media. In the case of repeated eigenvalues, the effective description of wave motion reveals three distinct asymptotic regimes depending on the symmetries of the germane eigenfunction basis. One of these regimes is shown to describe the so-called Dirac points, i.e. conical contacts between dispersion surfaces, that are relevant to the phenomenon of topological insulation. In situations involving nearby eigenvalues (some of which may be repeated), the leading-order model is found to invariably entail a \emph{Dirac-like} system of equations that generates ``blunted'' conical dispersion surfaces. We illustrate the analytical developments by several examples, including the Green's function near the edge of a band gap, and clusters of closely spaced dispersion branches. The mathematical framework established in this work naturally lends itself to the effective analysis of related (e.g. Maxwell or Navier) field equations and emerging physical phenomena, such as topologically protected states and topological networks~\cite{TWZ17,Makw18}, that revolve around the gapping of Dirac cones.

{\bf Funding.} This study was partially supported through the endowed Shimizu Professorship. 

{\bf Acknowledgments.} SM was partially supported by the NSF grant DMS-1439786 while in residence at ICERM Brown University during Fall 2017, and the AFOSR award FA9550-18-1-0131.


\end{document}